\documentclass[a4paper,11pt,reqno]{article}

\usepackage{amsmath,amsthm,amssymb}
\usepackage{thmtools}
\usepackage{relsize}
\usepackage[nospace,noadjust]{cite}
\usepackage[dvipsnames]{xcolor}
\usepackage[normalem]{ulem}
\usepackage[smalltableaux]{ytableau}
\usepackage{authblk}
\usepackage{enumitem}
\usepackage[margin=2.5cm]{geometry}
\usepackage{caption}
\usepackage{array}
\usepackage{stmaryrd}
\usepackage{scrextend}

\usepackage{tikz-cd}
\usepackage{tikz-3dplot}
\usetikzlibrary{calc,fit,matrix,arrows,automata,positioning,shadows}
\usetikzlibrary{decorations.pathreplacing,angles,quotes}

\definecolor{fast-rub-blau}{RGB}{0,70,146}
\usepackage{hyperref}
\hypersetup{
  colorlinks = true,
  urlcolor   = blue,
  linkcolor  = fast-rub-blau,
  citecolor  = fast-rub-blau
}
\usepackage[capitalize]{cleveref}

\declaretheorem[style=definition, numberwithin=section]{theorem}
\declaretheorem[style=definition, sibling=theorem]{definition, example,
  notation, remark, conjecture, corollary, proposition, lemma, problem}

\numberwithin{equation}{section}
\allowdisplaybreaks

\newcolumntype{x}[1]{>{\centering\arraybackslash\hspace{0pt}}p{#1}}

\newcommand{\mC}{\mathcal C}

\newcommand{\mM}{\mathcal{M}}
\newcommand{\mL}{\mathcal{L}}
\newcommand{\mU}{\mathcal{U}}
\newcommand{\mF}{\mathcal{F}}

\newcommand{\mI}{\mathcal{I}}
\newcommand{\mZ}{\mathcal{Z}}
\newcommand{\mS}{\mathcal{S}}
\newcommand{\mT}{\mathcal{T}}
\newcommand{\mN}{\mathcal{N}}

\newcommand{\F}{\mathbb F}
\newcommand{\R}{\mathbb R}

\newcommand{\mB}{\mathcal{B}}

\newcommand{\YL}{\mathfrak{Y}}

\newcommand{\ZZ}{\mathbb{Z}}

\DeclareMathOperator{\supp}{supp}

\DeclareMathOperator{\cyc}{cyc}

\DeclareMathOperator{\cl}{cl}

\DeclareMathOperator{\colsp}{colsp}

\DeclareMathOperator{\rank}{rank}

\newcommand{\Ht}{\operatorname{ht}}

\newcommand{\down}{{\downarrow}}

\title{Latroids on Young's Lattice}

\author[1]{Gianira N. Alfarano}
\author[2]{Sebastian Degen}

\affil[1]{Universit\'e de Rennes, IRMAR}
\affil[2]{Universit\"{a}t Bielefeld, Fakult\"{a}t f\"{u}r Mathematik, Bielefeld, Germany}

\affil[ ]{{\small \texttt{gianira-nicoletta.alfarano@univ-rennes.fr},
\texttt{sdegen@math.uni-bielefeld.de}}}

\date{}

\begin{document}

\maketitle
\begin{abstract}
Latroids are a generalization of matroids and polymatroids, given by a submodular rank function on a modular lattice. We study integer-valued latroids whose length function is the height of a finite modular lattice, without assuming that the lattice is complemented. We show that their cyclic flats form a lattice and, together with their ranks, determine the rank function. We also characterize which families of lattice elements, together with prescribed rank values, occur as the cyclic flats of such latroids. We then focus on principal intervals of Young's lattice, where we classify the latroids with modular rank functions and characterize the sparse paving ones. Finally, we attach a latroid to every Ferrers diagram rank-metric code and show that its cyclic flats, together with their ranks, determine the dimensions of all shortenings to Young subdiagrams.
\end{abstract}

\vspace{1em}


\tableofcontents

\medskip

\section{Introduction}
Latroids were introduced by Vertigan in \cite{vertigan2004latroids} to study
codes over modules within a generalization of matroids and
polymatroids. Their definition consists of a submodular rank function on a modular lattice, together with a strictly increasing modular length
function.
In \cite{gorla2025latroids}, Gorla and Salizzoni develop this framework further, relating latroids to code invariants and establishing cryptomorphic descriptions in terms of independent elements, bases, circuits and flats when the ambient lattice is complemented and the length is its height function.

Young's lattice $\YL$ is the set of all integer partitions, identified with their Young diagrams and ordered by inclusion; see \cite{stanley2012enumerative}. Its finite intervals are distributive lattices that are, in general, not complemented. Each interval $[Q,P]$ is isomorphic to the lattice of order ideals of the boxes of the skew diagram $P\setminus Q$, ordered coordinatewise. Products of chains, including those underlying the chain-support latroids of \cite[Section~3]{gorla2025latroids}, can also be described as intervals of Young's lattice. In this paper, we work with the principal intervals $\YL(P)=[\emptyset,P]$. Each of these intervals has a single atom whenever $P\neq\emptyset$, and is not complemented when $|P|\geq2$.

On finite distributive lattices, the integer-valued latroids whose length is the lattice height are
the \emph{Faigle-WSRFs} (where WSRF means ``weakly submodular rank function'') defined by Wild in~\cite{wild2008weakly}, arising from Faigle's work on pregeometries on partially ordered sets; see \cite{faigle1980geometries}.
Faigle's pregeometries also admit descriptions using notions of independence different from the classical condition $\rho(I)=\Ht(I)$.
For normalized, monotone, submodular integer-valued functions on a finite lattice whose values increase by at most one along
each cover, the flats form a graded upper semimodular lattice and determine the rank function; see
\cite{faigle1980geometries,wild2005rank}.

Cyclic flats provide another way of describing rank functions. Their lattice structure and properties have been established for matroids in~\cite{bonin2008lattice}, polymatroids in~\cite{csirmaz2020cyclic}, and $q$-matroids in~\cite{alfarano2024cyclic,gluesing2024decompositions}. Moreover, in \cite{byrne2025cyclic} this theory has been extended to $\mathcal{L}$-polymatroids on finite complemented modular lattices, where cyclic-flat
ranks together with atomic ranks determine the entire rank function.

In this paper, we study integer-valued latroids whose length function is the height of a finite modular lattice, and we do not assume that the lattice is complemented. Our main results are about their cyclic flats. We show that the cyclic flats form a lattice in \Cref{thm:cyclic-flats-lattice} and that, together with their ranks, they determine the rank function in \Cref{thm:reconstruction}. We also characterize the families of cyclic flats with ranks that arise in this way in \Cref{thm:cyclic-flat-axioms}. This gives a cryptomorphic description of these latroids by their cyclic flats. If the lattice is complemented, these latroids are $\mathcal{L}$-polymatroids, and results of this kind are proved in \cite[Sections~4 and~5]{byrne2025cyclic}. If the lattice is distributive, we show that every such latroid is the restriction of a matroid rank function to the order ideals of a poset, and that its cyclic flats are exactly the cyclic flats of a canonical matroid as such; see \Cref{prop:convolution}. Hence, for cyclic flats, we get new results for the modular lattices that are neither distributive nor complemented, such as the lattice of \Cref{ex:reconstruction-M3} and the submodule lattices of \Cref{rem:chain-rings}.

We then study these latroids on principal intervals $\YL(P_n)$ of Young's lattice, where $P_n$ is a partition of the integer $n$. Here, independent elements, circuits and bases do not need to determine the rank function, and a latroid of positive rank does not need to have a basis; see \Cref{ex:chain-nongeometric} and \Cref{rem:no-basis}. We classify the modular rank functions on $\YL(P_n)$. They count the boxes of a partition that lie in a fixed set $W$, and their independent partitions, bases, circuits, flats and cyclic flats can be read from $W$; see \Cref{prop:modular-classification,prop:modular-structure}. The size of the \emph{maximal hook} and the \emph{Durfee rank} of a partition are two examples of such modular rank functions. We also construct paving latroids on $\YL(P_n)$, and we show that this construction gives exactly the sparse paving ones; see \Cref{thm:paving-characterisation}. These are exactly the paving latroids whose dual is paving; see \Cref{prop:sparse-paving-dual}.

As an application to coding theory, we attach a latroid to every Ferrers diagram rank-metric code $C$, which is a space of matrices supported on a Ferrers diagram. Its rank function assigns to a Young subdiagram $P$ the value $|P|-\dim(C(P))$, where $C(P)$ is the shortening of $C$ to $P$. We describe its independent elements, circuits, flats and cyclic elements in terms of these shortenings, and we deduce from \Cref{thm:reconstruction} that the cyclic flats, together with their ranks, determine the dimensions of all shortenings to Young subdiagrams; see \Cref{thm:shortening-profile}. We also show that not every latroid on a principal interval of Young's lattice comes from a code; see \Cref{prop:not-from-code}. Finally, we look at the cyclic flats of a latroid as Ferrers diagrams, and we study the $q$-polymatroids defined by the corresponding spaces of matrices. Inclusions between diagrams imply a rank-difference inequality between the associated $q$-polymatroids, and when these are $q$-matroids, it says that one is a $q$-quotient of the other; see \Cref{prop:rank-difference}. This last property relates our setting also to the $q$-analogues of $\Delta$-matroids studied in \cite{CJ25qanadeltamats}.

The paper is organized as follows. In \Cref{sec:prelim}, we recall some notions on lattices, Young's lattice and latroids. In \Cref{sec:flats-cyclic}, we study flats and cyclic flats without complementedness. In \Cref{sec:young}, we study $\ZZ$-latroids on principal intervals of Young's lattice. \Cref{sec:ferrers} is about Ferrers diagram rank-metric codes, and \Cref{sec:qpoly} is about $q$-polymatroids. In Appendix~\ref{sec:enumeration} we enumerate all the $\ZZ$-latroids for $n\leq 4$.

\section{Preliminaries}\label{sec:prelim}

In this section, we collect some well-known facts about lattices, we introduce Young's lattice and we recall the basic notions on latroids.

\subsection{Lattices}

In this subsection, we recall the basic facts and definitions on lattices that we will need throughout the paper. We refer the reader to \cite{birkhoff1967lattice, stanley2012enumerative} for further details.

A \textbf{lattice} is a partially ordered set (poset) in which any two elements $X,Y$ have a least upper bound $X\vee Y$ (their \textbf{join}) and a greatest lower bound $X\wedge Y$ (their \textbf{meet}).

\begin{definition}\label{def:lattices}
    Let $\mL$ be a finite lattice with meet $\wedge$ and join $\vee$, least element $0_\mL$ and greatest element $1_\mL$. Let $X,Y\in\mL$.
    \begin{enumerate}
        \item[(1)] Let $X\leq Y$. The \textbf{interval} $[X,Y]$ is the set of all $Z\in\mL$ such that $X\leq Z\leq Y$.
        \item[(2)] We say that $Y$ is a \textbf{cover} of $X$, and we write $X\lessdot Y$, if $X<Y$ and there is no $Z\in\mL$ with $X<Z<Y$. An \textbf{atom} of $\mL$ is a cover of $0_\mL$.
        \item[(3)] $\mL$ is called \textbf{modular} if $X\vee(Y\wedge Z)=(X\vee Y)\wedge Z$ for all $X,Y,Z\in\mL$ with $X\leq Z$. It is called \textbf{distributive} if $X\wedge(Y\vee Z)=(X\wedge Y)\vee(X\wedge Z)$ for all $X,Y,Z\in\mL$.
        \item[(4)] $\mL$ is called \textbf{complemented} if for every $X\in\mL$ there exists $Y\in\mL$ with $X\wedge Y=0_\mL$ and $X\vee Y=1_\mL$.
        \item[(5)] $\mL$ is called \textbf{upper semimodular} if $X\wedge Y\lessdot X$ implies $Y\lessdot X\vee Y$, and \textbf{atomistic} if every element of $\mL$ is a join of atoms. It is called \textbf{geometric} if it is upper semimodular and atomistic.
    \end{enumerate}
\end{definition}

Every distributive lattice is modular. If $\mL$ is modular, then it satisfies the Jordan--Dedekind chain condition, that is, all maximal chains between two given elements of $\mL$ have the same length. Hence $\mL$ is \textbf{graded}, and the \textbf{height} $\Ht(X)$ of $X\in\mL$, which is the length of a maximal chain from $0_\mL$ to $X$, is well defined. The height function satisfies $\Ht(0_\mL)=0$, it is strictly increasing, and $\Ht(Y)=\Ht(X)+1$ whenever $X\lessdot Y$. Moreover, it is modular, that is, $\Ht(X)+\Ht(Y)=\Ht(X\vee Y)+\Ht(X\wedge Y)$ for all $X,Y\in\mL$.

\begin{notation}
    Throughout this paper, unless stated otherwise, $\mL$ denotes a finite modular lattice with least element $0_\mL$, greatest element $1_\mL$ and height function $\Ht$.
\end{notation}

\subsection{Young's lattice}\label{subsec:young-lattice}

In this subsection, we introduce Young's lattice, whose principal intervals are the main examples in this paper.

\begin{definition}
    A \textbf{partition} is a finite sequence $P=(\lambda_1,\ldots,\lambda_t)$ of integers with $\lambda_1\geq\lambda_2\geq\cdots\geq\lambda_t>0$. The \textbf{size} of $P$ is $|P|:=\lambda_1+\cdots+\lambda_t$. If $|P|=n$, we write $P\vdash n$. The \textbf{empty partition} $\emptyset$ has size $0$.
\end{definition}

We write repeated parts as exponents.~For instance, $(3,1^2)=(3,1,1)$, and $(a^s)$ is the partition with $s$ parts equal to $a$. We identify a partition $P=(\lambda_1,\ldots,\lambda_t)$ with its \textbf{Young diagram}, which is the set of boxes $P=\{(i,j)\mid 1\leq i\leq t,\ 1\leq j\leq\lambda_i\}$. The box $(i,j)$ lies in row $i$ and column $j$, and the diagram is drawn with its rows aligned on the left and its first row on top. Thus $(a^s)$ is a rectangle with $s$ rows and $a$ columns.

Throughout the paper, we use the same notation for a partition and its Young diagram. In particular, $b\in P$ means that $b$ is a box of the diagram of $P$, and $Q\subseteq P$ means that the diagram of $Q$ is contained in that of $P$.

\begin{example}\label{ex:young-diagram}
    The partition $P=(4,2,1)$ has size $|P|=4+2+1=7$, so $P\vdash 7$.
    Its Young diagram is
    \begin{center}
        \ydiagram{4,2,1}
    \end{center}
    It consists of three rows containing $4$, $2$ and $1$
    boxes, respectively. For instance, the last box of the first row is $(1,4)$.
\end{example}

\begin{definition}
\textbf{Young's lattice} $\YL$ is the set of all partitions, ordered by inclusion of their Young diagrams. For a partition $P$, we write $\YL(P):=[\emptyset,P]$ for the interval of all partitions whose diagram is contained in $P$, and we call its elements the \textbf{Young subdiagrams} of $P$. We call $\YL(P)$ a \textbf{principal interval} of Young's lattice.
\end{definition}

\begin{example}\label{ex:young-32}
    Let $P=(3,2)$. The lattice $\YL(P)$ has nine elements, and its (Hasse) diagram is shown in \Cref{fig:young-32}. Each edge corresponds to adding one box.
    
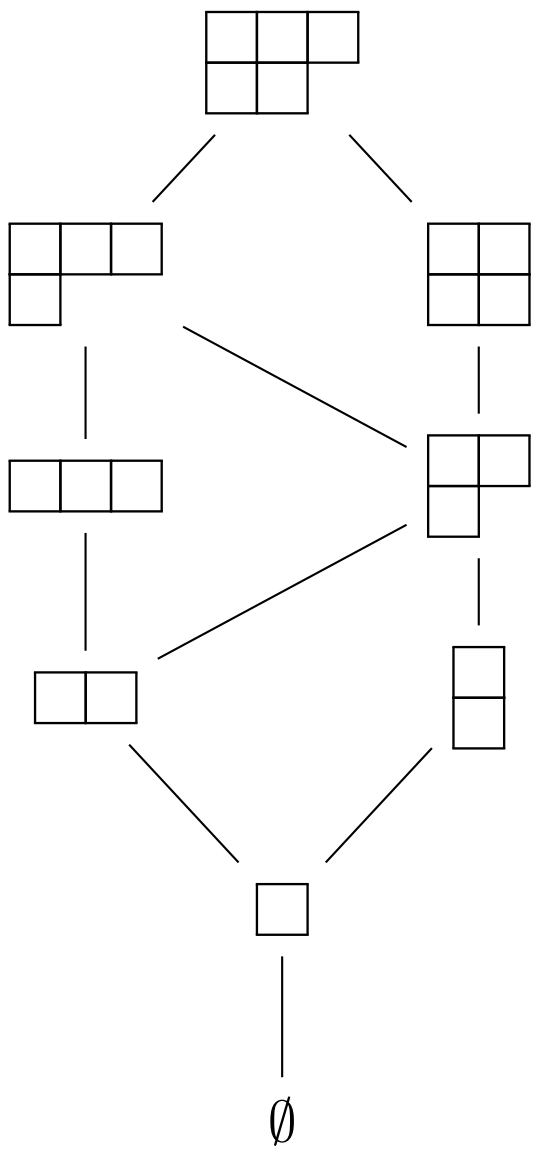
\begin{center}
        \begin{tikzpicture}[xscale=1.3,yscale=1.4]
            \node (e) at (0,0) {$\emptyset$};
            \node (1) at (0,1) {\ydiagram{1}};
            \node (2) at (-1,2) {\ydiagram{2}};
            \node (11) at (1,2) {\ydiagram{1,1}};
            \node (3) at (-1,3) {\ydiagram{3}};
            \node (21) at (1,3) {\ydiagram{2,1}};
            \node (31) at (-1,4) {\ydiagram{3,1}};
            \node (22) at (1,4) {\ydiagram{2,2}};
            \node (32) at (0,5) {\ydiagram{3,2}};
            \draw (e)--(1) (1)--(2) (1)--(11) (2)--(3) (2)--(21) (11)--(21) (3)--(31) (21)--(31) (21)--(22) (31)--(32) (22)--(32);
        \end{tikzpicture}
    \captionof{figure}{The Hasse diagram of the interval $\YL((3,2))$ of Young's lattice, where each partition is drawn as its Young diagram.}
        \label{fig:young-32}
    \end{center}
\end{example}

We regard boxes as elements of $\ZZ_{>0}^2$, ordered by $(i,j)\preceq(i',j')$ if $i\leq i'$ and $j\leq j'$. A finite set of boxes is a Young diagram if and only if it is an \textbf{order ideal} for $\preceq$, that is, together with any box $b$ it contains all boxes $b'\preceq b$. Thus the Young subdiagrams of $P$ are exactly the order ideals of the poset of boxes of $P$. For a finite set $S$ of boxes, we write
$$
    \down S:=\{b'\in\ZZ_{>0}^2\mid b'\preceq b\text{ for some }b\in S\}
$$
for its \textbf{downward closure}, which is the smallest Young diagram containing $S$. For a single box $b=(i,j)$, we let $R_b:=\down\{b\}$. This is the rectangle $(j^i)$ with $i$ rows and $j$ columns, and a Young diagram contains $b$ if and only if it contains $R_b$.

\begin{example}
    Let $S=\{(1,4),(3,2)\}$. Its downward closure is
    $\down S=R_{(1,4)}\cup R_{(3,2)}=(4)\cup(2^3)=(4,2,2)$.
    In the following diagram, the boxes of $S$ are marked
    with a dot, and all boxes of $\down S$ are shaded.
    
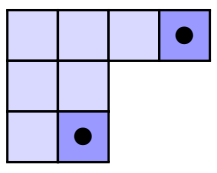
\begin{center}
        \begin{ytableau}
            *(blue!15) & *(blue!15) & *(blue!15) & *(blue!40) \bullet \\
            *(blue!15) & *(blue!15) \\
            *(blue!15) & *(blue!40) \bullet
        \end{ytableau}
    \end{center}
    Every Young diagram containing the two marked boxes
    must contain all the shaded boxes, so $(4,2,2)$ is the
    smallest Young diagram containing $S$.
\end{example}

The meet and the join in $\YL$ are the intersection and the union of diagrams. For instance,





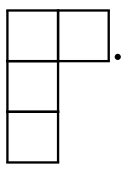
\[
    \ydiagram{1,1,1}\;\wedge\;\ydiagram{2,1}\; =\; \ydiagram{1,1}\quad\text{and}\quad
    \ydiagram{1,1,1}\;\vee\;\ydiagram{2,1}\; =\; \ydiagram{2,1,1}.
\]
Hence $\YL$ is an infinite distributive lattice, and so it is modular; see e.g.~\cite{birkhoff1967lattice}. Its least element is~$\emptyset$, and the height of $P$ is its size $|P|$, that is, the number of boxes of $P$.

\begin{remark}
Each interval $\YL(P)$ is a finite distributive lattice. If $P\neq\emptyset$, then $\YL(P)$ has exactly one atom, namely the one-box diagram $(1)$, because every nonempty diagram contains the box $(1,1)$. In particular, $\YL(P)$ is not complemented if $|P|\geq2$. Indeed, a complement $Y$ of $(1)$ would satisfy $Y\wedge(1)=\emptyset$, hence $Y=\emptyset$, and then $Y\vee(1)=(1)\neq P$.
\end{remark}

\begin{definition}
    Let $P=(\lambda_1,\ldots,\lambda_t)$ be a partition, and let $1\leq i\leq t$ and $1\leq j\leq\lambda_i$. The \textbf{hook} of $P$ at the box $(i,j)$ consists of the boxes of $P$ weakly to the right of $(i,j)$ in row $i$, together with the boxes of $P$ weakly below $(i,j)$ in column $j$, i.e.
    \[
        H_{(i,j)}^P
        :=
        \{(i,j')\in P\mid j'\geq j\}
        \cup
        \{(i',j)\in P\mid i'\geq i\}.
    \]
\end{definition}

For a nonempty partition $P=(\lambda_1,\ldots,\lambda_t)$, the \textbf{maximal hook} $H^P_{(1,1)}$ consists of the first row and the first column of $P$. So $|H^P_{(1,1)}|=\lambda_1+t-1$. The \textbf{diagonal hooks} of $P$ are the hooks~$H^P_{(i,i)}$ at the \textbf{diagonal boxes} $(i,i)\in P$.

\begin{example}\label{ex:hooks}
    Let $P=(3)$ and $P'=(2,1)$. Both are equal to their maximal hooks, and so are
    $$
        P\wedge P' = \ydiagram{2}\quad\text{and}\quad P\vee P' = \ydiagram{3,1}.
    $$
    Hence $|H^{P\wedge P'}_{(1,1)}|+|H^{P\vee P'}_{(1,1)}|=2+4=3+3=|H^{P}_{(1,1)}|+|H^{P'}_{(1,1)}|$.
\end{example}

\subsection{Latroids}

Latroids were introduced by Vertigan in \cite{vertigan2004latroids} as a generalization of matroids and polymatroids, motivated by the representation of codes over modules.~We use the notation of Gorla and Salizzoni \cite{gorla2025latroids}, who study latroids from the point of view of code invariants.

\begin{definition}[{\cite{vertigan2004latroids,gorla2025latroids}}]
    Let $A$ be a totally ordered abelian group. An \textbf{$A$-latroid on $\mL$} is a triple $L=(\rho,\lVert\cdot\rVert,\mL)$, consisting of a \textbf{rank function} $\rho\colon\mL\to A$ and a \textbf{length function} $\lVert\cdot\rVert\colon\mL\to A$, which satisfies the following axioms.
    \begin{enumerate}
        \item[($L1$)] $\rho(0_{\mL})=\lVert0_{\mL}\rVert=0_A$.
        \item[($L2$)] $\lVert\cdot\rVert$ is strictly increasing, that is, $\lVert X\rVert<\lVert Y\rVert$ for all $X,Y\in\mL$ with $X<Y$.
        \item[($L3$)] $\lVert\cdot\rVert$ is modular, that is, $\lVert X\rVert+\lVert Y\rVert=\lVert X\vee Y\rVert +\lVert X\wedge Y\rVert$ for all $X,Y\in\mL$.
        \item[($L4$)] $\rho$ is bounded increasing, that is, $0\leq\rho(Y)-\rho(X)\leq \lVert Y\rVert-\lVert X\rVert$ for all $X,Y\in\mL$ with $X<Y$.
        \item[($L5$)] $\rho$ is submodular, that is, $\rho(X)+\rho(Y)\geq\rho(X\vee Y)+\rho(X\wedge Y)$ for all $X,Y\in\mL$.
    \end{enumerate}
\end{definition}

By the properties of the height function, $(\rho,\Ht,\mL)$ is a $\ZZ$-latroid for every function $\rho\colon\mL\to\ZZ$ that satisfies $\rho(0_\mL)=0$, $(L4)$ and $(L5)$ with respect to $\Ht$.

\begin{definition}[{\cite{gorla2025latroids}}]
Let $L=(\rho,\lVert\cdot\rVert,\mL)$ be an $A$-latroid.
The \textbf{rank} of $L$ is $\rho(L):=\rho(1_\mL)$.
An element $I\in\mL$ is \textbf{independent} if
$\rho(I)=\lVert I\rVert$, and \textbf{dependent} otherwise.
A \textbf{basis} of $L$ is an element $B\in\mL$ such that
$\rho(B)=\lVert B\rVert=\rho(L)$.
A \textbf{circuit} of $L$ is a dependent element $C\in\mL$
such that every $C'<C$ is independent.
The \textbf{closure} of $F\in\mL$ is
$$
\cl(F):=\bigvee\{K\in\mL\mid \rho(F\vee K)=\rho(F)\}.
$$
An element $F\in\mL$ is a \textbf{flat} of $L$ if
$\cl(F)=F$. We write $\mF_L$ for the collection of flats of $L$.
\end{definition}

\begin{remark}\label{rem:no-basis}
    For a general latroid, a basis does not need to be the same as a maximal independent element, since a maximal independent element may have length strictly smaller than $\rho(1_{\mL})$. In particular, a latroid does not need to have a basis, although maximal independent elements always exist because $\mL$ is finite. For instance, on the chain $\YL((2))=\{\emptyset<(1)<(2)\}$ with length $|\cdot|$, the rank function with values $0,0,1$ defines a $\ZZ$-latroid of rank $1$. Its only independent element is $\emptyset$, so it has no basis.
\end{remark}

\begin{definition}\label{def:cyclic}
    Let $L=(\rho,\Ht,\mL)$ be a $\ZZ$-latroid. An element $C\in\mL$ is \textbf{cyclic} if $\rho(C')=\rho(C)$ for every $C'\lessdot C$. The least element $0_\mL$ is cyclic, since no element is covered by $0_\mL$. An element $C\in\mL$ is a \textbf{cyclic flat} if it is cyclic and a flat of $L$. We write $\mZ_L$ for the collection of cyclic flats of $L$. For $X\in\mL$, the \textbf{cyclic core} of $X$ is $\cyc(X):=\bigvee\{C\in\mL\mid C\leq X,\ C\text{ cyclic}\}$. We say that $L$ is \textbf{Schubert} if $\mZ_L$ is totally ordered.
\end{definition}

\begin{example}[{\cite[Example~2.3]{gorla2025latroids}}]
    Let $A$ be a totally ordered abelian group, and let $\lVert\cdot\rVert\colon\mL\to A$ be a strictly increasing modular function with $\lVert 0_\mL\rVert=0_A$. For $a\in A$ with $a\geq 0_A$, let $\rho_a(X):=\min\{a,\lVert X\rVert\}$. Then $(\rho_a,\lVert\cdot\rVert,\mL)$ is an $A$-latroid, which is called the \textbf{uniform $A$-latroid} and is denoted by $\mU_a(\mL,A)$.
\end{example}

\begin{remark}\label{rem:gs-complemented}
    For a $\ZZ$-latroid $(\rho,\Ht,\mL)$ on a complemented modular lattice $\mL$, Gorla and Salizzoni characterise the collection of flats $\mF_L$ by three conditions (F1)--(F3); see \cite[Proposition~2.30]{gorla2025latroids}, and (F1) and (F2) imply that $\mF_L$ is a lattice with the order induced by $\mL$. They also show that $\mF_L$ is a geometric lattice; see \cite[Theorem~2.31]{gorla2025latroids}. These results assume that $\mL$ is complemented, so they do not apply to the intervals $\YL(P)$ with $|P|\geq2$. Moreover, condition~(F3) is formulated in terms of atoms, and $\YL(P)$ has only one atom.
    In the next section, we recall what remains true for flats without the ``complemented'' assumption, and we extend the theory of cyclic flats to this setting.
\end{remark}

\section{Flats and cyclic flats on modular lattices}\label{sec:flats-cyclic}

In this section, we study $\ZZ$-latroids whose length function is the height. Throughout this section, $L=(\rho,\Ht,\mL)$ is such a $\ZZ$-latroid. Since $\Ht(Y)=\Ht(X)+1$ whenever $X\lessdot Y$, axiom $(L4)$ gives for all $X\lessdot Y$,
\begin{equation}\label{eq:unit-increase}
    \rho(Y)-\rho(X)\in\{0,1\}.
\end{equation}

In the terminology of Wild \cite{wild2005rank,wild2008weakly}, $\rho$ is a \emph{rank function} of $\mL$, that is, a monotone map with $\rho(0_\mL)=0$ and $\rho\leq\Ht$. Moreover, it is submodular and has \textbf{unit increase}, which means that~\eqref{eq:unit-increase} holds. If $\mL$ is distributive, these maps are exactly the \emph{Faigle-WSRFs} of \cite[p.~1002]{wild2008weakly}, and they are closely related to the pregeometries of Faigle; see \cite{faigle1980geometries}.
We will recall known results on flats. We extend the cyclic-flat results of Byrne and Fulcher \cite{byrne2025cyclic} to arbitrary finite modular lattices. The new cases are those in which the lattice is neither distributive nor complemented, since for distributive lattices the results also follow from the matroid case; see \Cref{prop:convolution}.

We need the following preliminary result.

\begin{lemma}\label{lem:cover-transfer}
    Let $\mL$ be a modular lattice, let $D\lessdot C$ in $\mL$, and let $X\in\mL$ with $X\leq C$ and $X\nleq D$. Then $X\vee D=C$ and $X\wedge D\lessdot X$.
\end{lemma}

\begin{proof}
Since $D\leq C$ and $X\leq C$, we have $D\leq X\vee D\leq C$. Moreover, $X\vee D\neq D$, since otherwise $X\leq D$. Hence $D<X\vee D\leq C$, and since $C$ covers $D$, we obtain $X\vee D=C$. For the second claim, note that $X\wedge D<X$, since $X\nleq D$. Let $Y\in\mL$ with $X\wedge D\leq Y\leq X$. By the first part, $D\leq Y\vee D\leq C$, so $Y\vee D\in\{D,C\}$. If $Y\vee D=D$, then $Y\leq X\wedge D$, and hence $Y=X\wedge D$. If $Y\vee D=C$, then the modular law, applied to $Y\leq X$, gives
$$ X = X\wedge C = X\wedge (Y\vee D) = Y\vee (X\wedge D) = Y.$$
Therefore $X\wedge D\lessdot X$.
\end{proof}

\subsection{Flats}\label{subsec:flats}

We recall the known facts about closure and flats due to Wild \cite{wild2005rank} for lattices of finite height and to Faigle \cite{faigle1980geometries} for distributive lattices. For complemented modular lattices they are also proved in \cite{gorla2025latroids}.

\begin{lemma}[{\cite[Theorem~1]{wild2005rank}}]\label{lem:closure}
    Let $L=(\rho,\Ht,\mL)$ be a $\ZZ$-latroid and let $X\in\mL$. The following hold.
    \begin{enumerate}
        \item[(1)] $\cl(X)$ is the largest element $Z\geq X$ with $\rho(Z)=\rho(X)$. In particular, $X\leq\cl(X)$ and $\rho(\cl(X))=\rho(X)$.
        \item[(2)] $\cl$ is a closure operator, and $\cl(X)$ is the smallest flat containing $X$, that is, $\cl(X)=\bigwedge\{F\in\mF_L\mid X\leq F\}$.
        \item[(3)] The following are equivalent.
        \begin{enumerate}
            \item[(i)] $X$ is a flat.
            \item[(ii)] $\rho(Z)>\rho(X)$ for every $Z>X$.
            \item[(iii)] $\rho(X')=\rho(X)+1$ for every $X'$ that covers $X$.
        \end{enumerate}
        \item[(4)] The meet of two flats is a flat. Consequently $\mF_L$, ordered by $\leq$, is a lattice with meet $F\wedge G$ and join $\cl(F\vee G)$. Its least element is $\cl(0_\mL)$ and its greatest element is $1_\mL$.
    \end{enumerate}
\end{lemma}

Wild works on lattices of finite height with weakly submodular maps; see \cite{wild2008weakly}. His closure $\bigvee\{Z\geq X\mid\rho(Z)=\rho(X)\}$ coincides with ours, and his flats are the elements described in (3)(ii). The equivalence of (3)(ii) and (3)(iii) follows from the monotonicity of $\rho$ and \eqref{eq:unit-increase}.~For distributive lattices see also \cite[Section~3]{faigle1980geometries}, and for complemented lattices see \cite[Propositions~2.29 and~2.30]{gorla2025latroids}.

\begin{theorem}[{\cite{faigle1980geometries,wild2005rank}}]\label{thm:flats-determine-rank}
    Let $L=(\rho,\Ht,\mL)$ be a $\ZZ$-latroid. The following hold.
    \begin{enumerate}
        \item[(1)] If $F,G\in\mF_L$ and $G$ covers $F$ in the lattice $\mF_L$, then $\rho(G)=\rho(F)+1$.
        \item[(2)] The lattice of flats $\mF_L$ is graded, with least element $\cl(0_\mL)$ and rank function $\rho|_{\mF_L}$, and it is upper semimodular.
        \item[(3)] For every $X\in\mL$ we have $\rho(X)=\Ht_{\mF_L}\bigl(\cl(X)\bigr)$, where $\Ht_{\mF_L}$ is the height function of $\mF_L$.
    \end{enumerate}
    In particular, two $\ZZ$-latroids $(\rho,\Ht,\mL)$ and $(\rho',\Ht,\mL)$ with the same flats coincide.
\end{theorem}

For lattices of finite height, (2) and (3) are due to Wild \cite[Theorems~2 and~3]{wild2005rank}. He shows that a weakly submodular map with unit increase is tight, which is the statement (3), and that its lattice of flats is graded and, by submodularity, upper semimodular. The upper semimodularity of the lattice of flats goes back to Gr\"atzer and Kiss \cite{gratzer1986construction}. Part (1) follows from (2). For distributive lattices these statements are in \cite{faigle1980geometries}, and for complemented lattices the argument is in the proof of \cite[Theorem~2.31]{gorla2025latroids}.

\begin{remark}\label{rem:flats-vs-complemented}
    If $\mL$ is complemented, the lattice of flats is geometric, and every finite geometric lattice occurs as the lattice of flats of such a $\ZZ$-latroid \cite[Theorem~2.31]{gorla2025latroids}. Without complementedness, the lattice of flats is still upper semimodular, but it does not need to be atomistic, and hence it does not need to be geometric, as \Cref{ex:chain-nongeometric} illustrates. In fact, every finite upper semimodular lattice is the lattice of flats of a $\ZZ$-latroid on a suitable distributive lattice; see \cite[Theorem~2, p.~31]{faigle1980geometries} and \cite[p.~1006]{wild2008weakly}. 
\end{remark}

\begin{example}\label{ex:chain-nongeometric}
    Let $\mL$ be the chain $0_\mL\lessdot X\lessdot Y\lessdot 1_\mL$, which is isomorphic to $\YL((3))$. Below we draw four functions $\rho$ on $\mL$, with their values in grey, and we highlight the flats by adding a box around them.
    $$
        \begin{tikzpicture}[yscale=0.9, inner sep=1.5pt,
            v/.style={label={[font=\scriptsize,gray]right:$#1$}},
            fl/.style={draw,rounded corners=2pt}]
            \node[v=0,fl] (a0) at (0,0) {$0_\mL$};
            \node[v=1] (a1) at (0,1) {$X$};
            \node[v=1,fl] (a2) at (0,2) {$Y$};
            \node[v=2,fl] (a3) at (0,3) {$1_\mL$};
            \draw (a0)--(a1)--(a2)--(a3);
            \node[font=\small] at (0,-1) {(a)};
            \node[v=0,fl] (b0) at (2.6,0) {$0_\mL$};
            \node[v=1] (b1) at (2.6,1) {$X$};
            \node[v=1] (b2) at (2.6,2) {$Y$};
            \node[v=1,fl] (b3) at (2.6,3) {$1_\mL$};
            \draw (b0)--(b1)--(b2)--(b3);
            \node[font=\small] at (2.6,-1) {(b)};
            \node[v=0] (c0) at (5.2,0) {$0_\mL$};
            \node[v=0] (c1) at (5.2,1) {$X$};
            \node[v=0,fl] (c2) at (5.2,2) {$Y$};
            \node[v=1,fl] (c3) at (5.2,3) {$1_\mL$};
            \draw (c0)--(c1)--(c2)--(c3);
            \node[font=\small] at (5.2,-1) {(c)};
            \node[v=0] (d0) at (7.8,0) {$0_\mL$};
            \node[v=0,fl] (d1) at (7.8,1) {$X$};
            \node[v=1] (d2) at (7.8,2) {$Y$};
            \node[v=1,fl] (d3) at (7.8,3) {$1_\mL$};
            \draw (d0)--(d1)--(d2)--(d3);
            \node[font=\small] at (7.8,-1) {(d)};
        \end{tikzpicture}
    $$
    Each of these functions increases by $0$ or $1$ along each cover. For $U<V$ in $\mL$,
$\rho(V)-\rho(U)$ is the sum of these increases over the $\Ht(V)-\Ht(U)$ covers
between $U$ and $V$, so $0\leq\rho(V)-\rho(U)\leq\Ht(V)-\Ht(U)$, which is $(L4)$. $(L5)$ holds trivially on a chain. So each of them defines a $\ZZ$-latroid $(\rho,\Ht,\mL)$. Consider first (a). By \Cref{lem:closure}\,(3), its flats are $0_\mL$, $Y$ and $1_\mL$, while $X$ is not a flat. The lattice of flats is the chain $0_\mL<Y<1_\mL$. It is graded and upper semimodular, but it is not atomistic, since $1_\mL$ is not a join of atoms. So it is not geometric. The $\ZZ$-latroid (b) has the same independent elements as (a), namely $0_\mL$ and $X$, and the same unique circuit $Y$, but its flats are $0_\mL$ and $1_\mL$. Similarly, the $\ZZ$-latroids (c) and (d) have the same unique independent element $0_\mL$, the same unique circuit $X$, and no basis, but different flats. Thus, unlike the case of complemented modular lattices with the height as length function \cite{gorla2025latroids},
the independent elements, circuits and bases do not determine a $\ZZ$-latroid in general.
\end{example}

\subsection{Cyclic flats}\label{subsec:cyclic-flats}

We now study the cyclic flats of a $\ZZ$-latroid $L=(\rho,\Ht,\mL)$ on a finite modular lattice $\mL$. We show that they form a lattice and that, together with their ranks, they determine the rank function. We first collect some properties of cyclic elements.

\begin{lemma}\label{lem:cyclic}
    Let $L=(\rho,\Ht,\mL)$ be a $\ZZ$-latroid. The following hold.
    \begin{enumerate}
        \item[(1)] If $C$ is cyclic, $C\leq X$, and $H\lessdot X$ with $C\nleq H$, then $\rho(H)=\rho(X)$.
        \item[(2)] If $C_1,C_2$ are cyclic, then $C_1\vee C_2$ is cyclic. Consequently, $\cyc(X)$ is the largest cyclic element below $X$.
        \item[(3)] If $C$ is cyclic, then $\cl(C)$ is cyclic.
        \item[(4)] If $F$ is a flat, then $\cyc(F)$ is a flat, hence a cyclic flat.
    \end{enumerate}
\end{lemma}

\begin{proof}
    (1) By \Cref{lem:cover-transfer}, applied to $H\lessdot X$ and $C$, we have $C\vee H=X$ and $C\wedge H\lessdot C$. So $\rho(C\wedge H)=\rho(C)$, because $C$ is cyclic. Submodularity gives $\rho(X)-\rho(H)\leq\rho(C)-\rho(C\wedge H)=0$. Since $H\leq X$, the monotonicity of $\rho$ gives the reverse
inequality, and hence  $\rho(H)=\rho(X)$.

    (2) Let $C:=C_1\vee C_2$ and let $H\lessdot C$. Since $H<C$, the elements $C_1$ and $C_2$ do not both lie below $H$. Say $C_1\nleq H$. Then (1), applied to $C_1\leq C$, gives $\rho(H)=\rho(C)$. So $C$ is cyclic. Since $0_\mL$ is cyclic and $\mL$ is finite, $\cyc(X)$ is a finite join of cyclic elements, and hence it is the largest cyclic element below $X$.

    (3) Let $F:=\cl(C)$ and let $H\lessdot F$. If $C\leq H$, then $\rho(F)=\rho(C)\leq\rho(H)\leq\rho(F)$ by \Cref{lem:closure}\,(1). Otherwise, (1) applied to $C\leq F$ gives $\rho(H)=\rho(F)$. Hence $F$ is cyclic.

    (4) Let $K:=\cyc(F)$. By (2) and (3), $\cl(K)$ is cyclic, and $\cl(K)\leq F$ by \Cref{lem:closure}\,(2). By the maximality of $K$ we get $\cl(K)\leq K$. Hence $K=\cl(K)$ is a flat.
\end{proof}

\begin{theorem}\label{thm:cyclic-flats-lattice}
    Let $L=(\rho,\Ht,\mL)$ be a $\ZZ$-latroid. Then $\mZ_L$, ordered by $\leq$, is a lattice. Its join is $Z_1\vee_{\mZ}Z_2=\cl(Z_1\vee Z_2)$ and its meet is $Z_1\wedge_{\mZ}Z_2=\cyc(Z_1\wedge Z_2)$. Its least element is $\cl(0_\mL)$ and its greatest element is $\cyc(1_\mL)$.
\end{theorem}

\begin{proof}
    Let $Z_1,Z_2\in\mZ_L$. By \Cref{lem:cyclic}\,(2) and (3), $\cl(Z_1\vee Z_2)$ is cyclic, and it is a flat by \Cref{lem:closure}\,(2). If $Z\in\mZ_L$ lies above $Z_1$ and $Z_2$, then $Z\geq Z_1\vee Z_2$, and hence $Z\geq\cl(Z_1\vee Z_2)$. Dually, $Z_1\wedge Z_2$ is a flat by \Cref{lem:closure}\,(4). So $\cyc(Z_1\wedge Z_2)$ is a cyclic flat by \Cref{lem:cyclic}\,(4), and every cyclic flat below $Z_1$ and $Z_2$ lies below $\cyc(Z_1\wedge Z_2)$ by \Cref{lem:cyclic}\,(2). Finally, $\cl(0_\mL)$ is cyclic by \Cref{lem:cyclic}\,(3) and lies below every flat, and $\cyc(1_\mL)$ is a cyclic flat by \Cref{lem:cyclic}\,(4) and lies above every cyclic element.
\end{proof}

For $X\in\mL$, we call $\nu(X):=\Ht(X)-\rho(X)$ the \textbf{nullity} of $X$.
We include a proof of the following identity for completeness.

\begin{lemma}\label{lem:nullity}
    Let $L=(\rho,\Ht,\mL)$ be a $\ZZ$-latroid. Then $\nu$ is increasing,
    and $\nu(X)=\nu(\cyc(X))$ for every $X\in\mL$.
\end{lemma}

\begin{proof}
    By $(L4)$, $\nu$ is increasing. By \eqref{eq:unit-increase}, we have $\nu(Y)-\nu(D)\in\{0,1\}$ for $D\lessdot Y$, with $\nu(Y)=\nu(D)$ if and only if $\rho(Y)=\rho(D)+1$. Let $X\in\mL$, and choose $W\leq X$ minimal with $\nu(W)=\nu(X)$. If $D\lessdot W$, then $\nu(D)\neq\nu(W)$ by minimality, so $\nu(D)=\nu(W)-1$, that is, $\rho(D)=\rho(W)$. Hence $W$ is cyclic. So $W\leq\cyc(X)$, and $\nu(X)=\nu(W)\leq\nu(\cyc(X))\leq\nu(X)$.
\end{proof}

\begin{theorem}\label{thm:reconstruction}
    Let $L=(\rho,\Ht,\mL)$ be a $\ZZ$-latroid. Then for every $X\in\mL$,
    \[
        \rho(X)=\min_{Z\in\mZ_L}\Bigl(\rho(Z)+\Ht(X\vee Z)-\Ht(Z)\Bigr),
    \]
    and the minimum is attained at $Z=\cyc(\cl(X))$. In particular, $L$ is determined by its cyclic flats together with their ranks.
\end{theorem}

\begin{proof}
    For every $Z\in\mZ_L$, the monotonicity of $\rho$ and $(L4)$ give $\rho(X)\leq\rho(X\vee Z)\leq\rho(Z)+\Ht(X\vee Z)-\Ht(Z)$.
    For the converse inequality, let $\nu$ be the nullity, and let $F:=\cl(X)$ and $Z:=\cyc(F)$. By \Cref{lem:closure}\,(2) and \Cref{lem:cyclic}\,(4), $Z$ is a cyclic flat. Since $Z\leq X\vee Z\leq F$, \Cref{lem:nullity} gives $\nu(Z)\leq\nu(X\vee Z)\leq\nu(F)=\nu(\cyc(F))=\nu(Z)$. Therefore $\rho(X\vee Z)=\rho(Z)+\Ht(X\vee Z)-\Ht(Z)$. Finally, $X\leq X\vee Z\leq\cl(X)$ gives $\rho(X\vee Z)=\rho(X)$ by \Cref{lem:closure}\,(1).
\end{proof}

Together, \Cref{thm:cyclic-flats-lattice,thm:reconstruction} lead to a description of $\ZZ$-latroids by their cyclic flats together with their ranks, which is analogous to the axioms for the cyclic flats of matroids from \cite{bonin2008lattice} and of $q$-matroids from \cite{alfarano2024cyclic}.

\begin{theorem}\label{thm:cyclic-flat-axioms}
    Let $\mL$ be a finite modular lattice, let $\mZ\subseteq\mL$ be nonempty, and let $r\colon\mZ\to\ZZ$. Then there is a $\ZZ$-latroid $L=(\rho,\Ht,\mL)$ with $\mZ_L=\mZ$ and $\rho(Z)=r(Z)$ for every $Z\in\mZ$ if and only if the following hold.
    \begin{enumerate}
        \item[(Z0)] $\mZ$, with the order induced by $\mL$, is a lattice. We write $\vee_\mZ$ and $\wedge_\mZ$ for its join and meet, and $0_\mZ$ for its least element.
        \item[(Z1)] $r(0_\mZ)=0$.
        \item[(Z2)] $0<r(Z')-r(Z)<\Ht(Z')-\Ht(Z)$ for all $Z,Z'\in\mZ$ with $Z<Z'$.
        \item[(Z3)] $r(Z)+r(Z')\geq r(Z\vee_\mZ Z')+r(Z\wedge_\mZ Z')+\Ht(Z\wedge Z')-\Ht(Z\wedge_\mZ Z')$ for all $Z,Z'\in\mZ$.
    \end{enumerate}
    In this case, $L$ is unique, and its rank function is
    $$
        \rho(X)=\min_{Z\in\mZ}\Bigl(r(Z)+\Ht(X\vee Z)-\Ht(Z)\Bigr).
    $$
\end{theorem}

\begin{proof}
    Suppose first that $L=(\rho,\Ht,\mL)$ is a $\ZZ$-latroid with $\mZ_L=\mZ$ and $\rho(Z)=r(Z)$ for every $Z\in\mZ$. By \Cref{thm:cyclic-flats-lattice}, $\mZ$ is a lattice with least element $\cl(0_\mL)$, so (Z0) holds, and (Z1) holds since $\rho(\cl(0_\mL))=\rho(0_\mL)=0$ by \Cref{lem:closure}\,(1). Let $Z<Z'$ in $\mZ$. Since $Z$ is a flat, $r(Z')>r(Z)$ by \Cref{lem:closure}\,(3). Choose $H\lessdot Z'$ with $Z\leq H$. Since $Z'$ is cyclic, $r(Z')=\rho(H)$, and $(L4)$ gives $r(Z')-r(Z)=\rho(H)-\rho(Z)\leq\Ht(H)-\Ht(Z)=\Ht(Z')-\Ht(Z)-1$. This proves (Z2). For (Z3), let $Z,Z'\in\mZ$, and let $C:=Z\vee_\mZ Z'$ and $D:=Z\wedge_\mZ Z'$. By \Cref{thm:cyclic-flats-lattice}, $C=\cl(Z\vee Z')$ and $D=\cyc(Z\wedge Z')$. Hence $\rho(Z\vee Z')=r(C)$ by \Cref{lem:closure}\,(1), and $\Ht(Z\wedge Z')-\rho(Z\wedge Z')=\Ht(D)-r(D)$ by \Cref{lem:nullity}. So (Z3) is the inequality $(L5)$ for $Z$ and $Z'$. Finally, $L$ is unique, and $\rho$ is given by the formula above, by \Cref{thm:reconstruction}.

        Conversely, suppose that (Z0)--(Z3) hold. For $Z\in\mZ$ and $X\in\mL$, let
    $$
        g_Z(X):=r(Z)+\Ht(X\vee Z)-\Ht(Z)=r(Z)+\Ht(X)-\Ht(X\wedge Z),
    $$
    where the second equality holds since $\Ht$ is modular. Let $\rho(X):=\min_{Z\in\mZ}g_Z(X)$. We show that $L=(\rho,\Ht,\mL)$ is a $\ZZ$-latroid with $\mZ_L=\mZ$ and $\rho(Z)=r(Z)$ for every $Z\in\mZ$.

    We first show that $L$ is a $\ZZ$-latroid. By (Z1) and (Z2), we have $r(Z)\geq0$ for every $Z\in\mZ$, with equality for $Z=0_\mZ$. Hence $\rho(0_\mL)=\min_{Z\in\mZ}r(Z)=0$. Let $X\lessdot X'$ and $Z\in\mZ$. Then $(X'\wedge Z)\wedge X=X\wedge Z$, and $X\leq(X'\wedge Z)\vee X\leq X'$. So, the modularity of $\Ht$ gives $\Ht(X'\wedge Z)-\Ht(X\wedge Z)=\Ht((X'\wedge Z)\vee X)-\Ht(X)\in\{0,1\}$. Hence $g_Z(X)\leq g_Z(X')\leq g_Z(X)+1$ for every $Z\in\mZ$, and taking the minimum over $Z$ gives $\rho(X)\leq\rho(X')\leq\rho(X)+1$. Now, let $X<Y$ in $\mL$, and choose a maximal chain
$X=X_0\lessdot X_1\lessdot\cdots\lessdot X_m=Y$. Since $\Ht$ increases by exactly $1$
along each cover, $m=\Ht(Y)-\Ht(X)$ and, for every $i$, $\rho(X_i)-\rho(X_{i-1})\in\{0,1\}$, so
$$
    0\leq\rho(Y)-\rho(X)=\sum_{i=1}^{m}\bigl(\rho(X_i)-\rho(X_{i-1})\bigr)\leq m=\Ht(Y)-\Ht(X).
$$
This proves $(L4)$. For $(L5)$, let $X,Y\in\mL$, and let $A,B\in\mZ$ with $\rho(X)=g_A(X)$ and $\rho(Y)=g_B(Y)$. Let $C:=A\vee_\mZ B$, $D:=A\wedge_\mZ B$ and $T:=X\wedge Y$. Since $A,B\leq C$, we have $(X\wedge A)\vee(Y\wedge B)\leq(X\vee Y)\wedge C$. Moreover, $(X\wedge A)\wedge(Y\wedge B)=T\wedge A\wedge B$. By modularity of $\Ht$,
    $$
        \Ht(X\wedge A)+\Ht(Y\wedge B)\leq\Ht((X\vee Y)\wedge C)+\Ht(T\wedge A\wedge B).
    $$
    Since $D\leq A\wedge B$, modularity also gives $\Ht(T\wedge A\wedge B)-\Ht(T\wedge D)=\Ht((T\wedge A\wedge B)\vee D)-\Ht(D)\leq\Ht(A\wedge B)-\Ht(D)$. Together with (Z3) for $A$ and $B$, and with $\Ht(X)+\Ht(Y)=\Ht(X\vee Y)+\Ht(T)$, we get
    $$
        \rho(X)+\rho(Y)=g_A(X)+g_B(Y)\geq g_C(X\vee Y)+g_D(X\wedge Y)\geq\rho(X\vee Y)+\rho(X\wedge Y).
    $$
    So $L$ is a $\ZZ$-latroid.

    Next, we show that $\rho(Z)=r(Z)$ for every $Z\in\mZ$. Let $Z,A\in\mZ$, and let $C:=Z\vee_\mZ A$ and $D:=Z\wedge_\mZ A$. By (Z3),
    $$
        g_A(Z)=r(A)+\Ht(Z)-\Ht(Z\wedge A)\geq r(C)+\bigl(\Ht(Z)-\Ht(D)-r(Z)+r(D)\bigr).
    $$
    Since $D\leq Z\leq C$, (Z2) gives $r(C)\geq r(Z)$, with strict inequality if $C\neq Z$. It also shows that the term in brackets is nonnegative, and positive if $D\neq Z$. If $A\neq Z$, then $C\neq Z$ or $D\neq Z$. So $g_A(Z)\geq r(Z)+1$ for every $A\neq Z$, while $g_Z(Z)=r(Z)$. In particular, $\rho(Z)=r(Z)$.

    We now show that every $Z\in\mZ$ is a cyclic flat of $L$. Recall that every $g_A$ changes by $0$ or $1$ along a cover. If $Z\lessdot X$, then $g_Z(X)=r(Z)+1$, and $g_A(X)\geq g_A(Z)\geq r(Z)+1$ for $A\neq Z$. So $\rho(X)=r(Z)+1>\rho(Z)$, and $Z$ is a flat by \Cref{lem:closure}\,(3). If $H\lessdot Z$, then $g_Z(H)=r(Z)$, and $g_A(H)\geq g_A(Z)-1\geq r(Z)$ for $A\neq Z$. So $\rho(H)=r(Z)=\rho(Z)$, and $Z$ is cyclic.

    Finally, we show that every cyclic flat of $L$ lies in $\mZ$. Let $F\in\mZ_L$, and choose $A\in\mZ$ with $\rho(F)=g_A(F)$. Suppose that $F\nleq A$. Then $F\wedge A<F$, and we can choose $H\lessdot F$ with $F\wedge A\leq H$. Then $H\wedge A=F\wedge A$, so $g_A(H)=g_A(F)-1$, and hence $\rho(H)<\rho(F)$. This contradicts the fact that $F$ is cyclic. So $F\leq A$, and $\rho(F)=g_A(F)=r(A)=\rho(A)$. Since $F$ is a flat, \Cref{lem:closure}\,(1) gives $A\leq\cl(F)=F$, and hence $F=A\in\mZ$.
\end{proof}

Note that $\mZ$ does not need to be a sublattice of $\mL$, since $\vee_\mZ$ and $\wedge_\mZ$ may differ from $\vee$ and~$\wedge$. This is why the term $\Ht(Z\wedge Z')-\Ht(Z\wedge_\mZ Z')\geq0$ appears in (Z3).

If $\mL$ is distributive, $\ZZ$-latroids with length the height are closely related to matroids. Recall that an \emph{order ideal} of a finite poset $B$ is a subset $I\subseteq B$ such that $x\in I$ and $y\leq x$ imply $y\in I$. The order ideals of $B$, ordered by inclusion, form a distributive lattice $J(B)$, whose join and meet are the union and the intersection, and whose height function is the size $|\cdot|$. It is known that every finite distributive lattice is of this form; see \cite{birkhoff1967lattice}. For instance, $\YL(P)$ is the lattice of order ideals of the set of boxes of $P$. For a matroid on $B$ with rank function $r$, recall that a set $F\subseteq B$ is a \textbf{flat} if $r(F\cup\{b\})>r(F)$ for every $b\in B\setminus F$, and that a set $Z\subseteq B$ is \textbf{cyclic} if $r(Z\setminus\{z\})=r(Z)$ for every $z\in Z$. These conditions involve all sets obtained by adding or removing one element, and not only the covers in $J(B)$.

\begin{proposition}\label{prop:convolution}
    Let $B$ be a finite poset and let $L=(\rho,|\cdot|,J(B))$ be a $\ZZ$-latroid. For $S\subseteq B$, let
    $$
        \widehat r(S):=\min_{I\in J(B)}\bigl(\rho(I)+|S\setminus I|\bigr).
    $$
    Then the following hold.
    \begin{enumerate}
        \item[(1)] $\widehat r$ is the rank function of a matroid $M$ on $B$, and $\widehat r(I)=\rho(I)$ for every $I\in J(B)$.
        \item[(2)] The cyclic flats of the matroid $M$ are exactly the cyclic flats of $L$. In particular, every cyclic flat of $M$ is an order ideal of $B$.
    \end{enumerate}
\end{proposition}

\begin{proof}
    (1) Let $I\in J(B)$. Choosing $I$ itself in the minimum gives $\widehat r(I)\leq\rho(I)$. Conversely, let $I'\in J(B)$. Then $I\cap I'\in J(B)$, and $(L4)$ gives $\rho(I)\leq\rho(I\cap I')+|I\setminus I'|\leq\rho(I')+|I\setminus I'|$. Hence $\widehat r(I)=\rho(I)$, and in particular $\widehat r(\emptyset)=0$. For $I\in J(B)$ and $S\subseteq B$, let $f_I(S):=\rho(I)+|S\setminus I|$. Then $f_I(S)\leq f_I(S\cup\{b\})\leq f_I(S)+1$ for every $b\in B$, and taking the minimum over $I$ gives $\widehat r(S)\leq\widehat r(S\cup\{b\})\leq\widehat r(S)+1$. It remains to show that $\widehat r$ is submodular. Let $S,T\subseteq B$, and choose $I,I'\in J(B)$ attaining $\widehat r(S)$ and $\widehat r(T)$, respectively. We claim that
    $$
        |(S\cup T)\setminus(I\cup I')|+|(S\cap T)\setminus(I\cap I')|\leq|S\setminus I|+|T\setminus I'|.
    $$
    Indeed, an element of $S\cap T$ contributes $2$, $1$ or $0$ to both sides, according to whether it lies in none, exactly one or both of $I$ and $I'$. An element of $S\setminus T$ contributes $1$ to the left-hand side only if it does not lie in $I\cup I'$, and then it also contributes $1$ to the right-hand side. The same holds for the elements of $T\setminus S$. Together with $(L5)$, the claim gives
    $$
        \widehat r(S\cup T)+\widehat r(S\cap T)\leq\rho(I\cup I')+\rho(I\cap I')+|(S\cup T)\setminus(I\cup I')|+|(S\cap T)\setminus(I\cap I')|\leq\widehat r(S)+\widehat r(T).
    $$

    (2) Let $\nu$ be the nullity of $L$, so $\nu(Y)=|Y|-\rho(Y)$ for $Y\in J(B)$. By \Cref{lem:nullity}, $\nu$ is increasing. Let $Z$ be a cyclic flat of $M$, and let $I\in J(B)$ attain $\widehat r(Z)$. If there were $z\in Z\setminus I$, then $\widehat r(Z\setminus\{z\})\leq\rho(I)+|Z\setminus I|-1=\widehat r(Z)-1$, which is impossible since $Z$ is cyclic in $M$. So $Z\subseteq I$, and $\widehat r(I)=\rho(I)=\widehat r(Z)$. Since $Z$ is a flat of $M$, we get $I\subseteq Z$, and hence $Z=I$ is an order ideal. By (1), every order ideal strictly containing $Z$ has larger rank than $Z$, since $Z$ is a flat of $M$, and $\rho(Y)=\rho(Z)$ for every $Y\lessdot Z$ in $J(B)$, since $Z$ is cyclic in $M$. By \Cref{lem:closure}\,(3) and \Cref{def:cyclic}, $Z$ is a cyclic flat of $L$.

    Conversely, let $Z$ be a cyclic flat of $L$. Suppose first that $\widehat r(Z\cup\{y\})=\rho(Z)$ for some $y\notin Z$, and let $I\in J(B)$ attain $\widehat r(Z\cup\{y\})$. By $(L5)$ and $(L4)$,
    $$
        \rho(Z\cup I)\leq\rho(I)+\rho(Z)-\rho(Z\cap I)\leq\rho(I)+|Z\setminus I|=\rho(Z)-|\{y\}\setminus I|.
    $$
    Since $\rho(Z\cup I)\geq\rho(Z)$, we get $y\in I$ and $\rho(Z\cup I)=\rho(Z)$. This contradicts the fact that $Z$ is a flat of $L$, since $Z\cup I$ is an order ideal strictly containing $Z$. So $Z$ is a flat of $M$. Suppose now that $\widehat r(Z\setminus\{z\})=\rho(Z)-1$ for some $z\in Z$, and let $I\in J(B)$ attain $\widehat r(Z\setminus\{z\})$. Let $K:=Z\cap I\in J(B)$. Since $\rho(K)\leq\rho(I)$, we get
    $$
        \nu(K)\geq|Z\cap I|-\rho(I)=|Z\cap I|-\rho(Z)+1+|(Z\setminus\{z\})\setminus I|.
    $$
    If $z\in I$, the right-hand side is $\nu(Z)+1$, which contradicts the monotonicity of $\nu$, since $K\subseteq Z$. If $z\notin I$, the right-hand side is $\nu(Z)$, and $K\subseteq Z\setminus\{z\}$. Since $Z$ is cyclic in $L$, every order ideal $Y$ strictly contained in $Z$ lies below some $Y'\lessdot Z$ with $\nu(Y')=\nu(Z)-1$, and hence $\nu(Y)<\nu(Z)$. This is again a contradiction. So $Z$ is cyclic in $M$.
\end{proof}

Hence, if $\mL$ is distributive, every $\ZZ$-latroid $(\rho,\Ht,\mL)$ is the restriction of a matroid rank function to the order ideals. Conversely, the restriction of any matroid rank function on $B$ to $J(B)$ defines a $\ZZ$-latroid, since the join and the meet are the union and the intersection. The function $\widehat r$ is a convolution; see \cite{csirmaz2020cyclic}. By \Cref{prop:convolution}\,(2), for distributive lattices \Cref{thm:cyclic-flats-lattice,thm:reconstruction} can also be deduced from the corresponding results for matroids \cite{bonin2008lattice}, applied to $\widehat r$. The same holds for \Cref{thm:cyclic-flat-axioms}. Indeed, for $\mL=J(B)$, conditions (Z0)--(Z3) are the axioms of \cite{bonin2008lattice}, and the rank function $S\mapsto\min_{Z\in\mZ}(r(Z)+|S\setminus Z|)$ of the matroid constructed there is the function $\widehat r$ of its restriction to $J(B)$. A modular lattice which is not distributive is not a lattice of order ideals, so this argument is not applicable in that case. We illustrate \Cref{thm:cyclic-flat-axioms,thm:reconstruction} on a lattice which is neither distributive nor complemented, so that neither the results of \cite{byrne2025cyclic} nor \Cref{prop:convolution} apply.

\begin{example}\label{ex:reconstruction-M3}
    Let $\mL$ be the lattice with elements $0_\mL,A,X_1,X_2,X_3,B,1_\mL$ and covers $0_\mL\lessdot A\lessdot X_i\lessdot B\lessdot 1_\mL$ for $i=1,2,3$. We draw it together with the values, in grey, of a function $\rho$, which we define below.
    $$
        
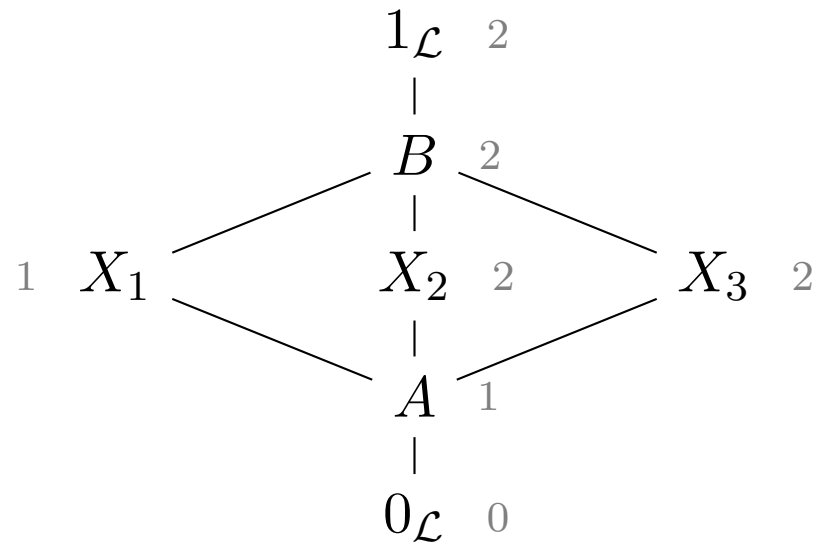
\begin{tikzpicture}[xscale=1.1,yscale=0.8]
            \node[label={[font=\scriptsize,gray]right:$0$}] (0) at (0,0) {$0_\mL$};
            \node[label={[font=\scriptsize,gray]right:$1$}] (A) at (0,1) {$A$};
            \node[label={[font=\scriptsize,gray]left:$1$}] (X1) at (-1.8,2) {$X_1$};
            \node[label={[font=\scriptsize,gray]right:$2$}] (X2) at (0,2) {$X_2$};
            \node[label={[font=\scriptsize,gray]right:$2$}] (X3) at (1.8,2) {$X_3$};
            \node[label={[font=\scriptsize,gray]right:$2$}] (B) at (0,3) {$B$};
            \node[label={[font=\scriptsize,gray]right:$2$}] (1) at (0,4) {$1_\mL$};
            \draw (0)--(A) (A)--(X1) (A)--(X2) (A)--(X3) (X1)--(B) (X2)--(B) (X3)--(B) (B)--(1);
        \end{tikzpicture}
    $$
    The interval $[A,B]$ is the diamond $M_3$ (see \cite{birkhoff1967lattice}), which is modular but not distributive. So $\mL$ is not distributive. It is modular, since it is obtained from $M_3$ by adjoining a new least element~$0_\mL$ and a new greatest element $1_\mL$, and the modular law holds trivially whenever one of the elements involved is $0_\mL$ or $1_\mL$. Moreover, $\mL$ is not complemented, since every element different from $0_\mL$ lies above $A$, and hence $A$ has no complement. The height function takes the values $0,1,2,2,2,3,4$ on $0_\mL,A,X_1,X_2,X_3,B,1_\mL$. We now use \Cref{thm:cyclic-flat-axioms} to construct a $\ZZ$-latroid on $\mL$. Let
    $\mZ:=\{0_\mL<X_1<1_\mL\}$, $r(0_\mL)=0$, $r(X_1)=1$, $r(1_\mL)=2$.
    Since $\mZ$ is a chain, (Z0) holds, and (Z3) holds with equality, because $Z\vee_\mZ Z'$ and $Z\wedge_\mZ Z'$ are the larger and the smaller of $Z$ and $Z'$. The heights of $0_\mL$, $X_1$ and $1_\mL$ are $0$, $2$ and $4$, so (Z1) and (Z2) hold as well. By \Cref{thm:cyclic-flat-axioms}, there is a unique $\ZZ$-latroid $L=(\rho,\Ht,\mL)$ with $\mZ_L=\mZ$ and $\rho(Z)=r(Z)$ for $Z\in\mZ$, and
    $$
        \rho(X)=\min\{\Ht(X),\;1+\Ht(X\vee X_1)-2,\;2\}.
    $$
    This gives the values $0,1,1,2,2,2,2$ on $0_\mL,A,X_1,X_2,X_3,B,1_\mL$, as in the picture. For instance, since $X_2\vee X_1=B$, we have $\rho(X_2)=\min\{2,\;1+3-2,\;2\}=2$, and $\rho(A)=\min\{1,\;1+2-2,\;2\}=1$.
    For $A$, the minimum is attained at $\cyc(\cl(A))=X_1$, as in \Cref{thm:reconstruction}, and also at $0_\mL$.
    \end{example}

\begin{remark}\label{rem:literature-cyclic}
    \Cref{lem:cyclic} and \Cref{thm:cyclic-flats-lattice} extend \cite[Lemmas~4.7, 4.8, 4.10, 4.13 and Proposition~4.14]{byrne2025cyclic} from complemented modular lattices to arbitrary finite modular lattices. Atoms are replaced by covers via \Cref{lem:cover-transfer}.~For integer-valued rank functions with unit increase, the cyclic elements of \cite[Definition~4.3]{byrne2025cyclic} are the cyclic elements of \Cref{def:cyclic}. In contrast, the flats of \cite[Definition~4.2]{byrne2025cyclic} are defined in terms of atoms. On $\YL(P)$, for $P\neq\emptyset$, the unique atom $(1)$ lies below every nonempty diagram, so every nonempty diagram would satisfy their definition of a flat.~\Cref{thm:reconstruction} is a counterpart of \cite[Theorem~4.25]{byrne2025cyclic} without complements. For matroids, $q$-matroids and polymatroids, see \cite{bonin2008lattice, csirmaz2020cyclic, alfarano2024cyclic}.

    The assumption that $\rho$ takes integer values is not needed for \Cref{lem:cyclic} and \Cref{thm:cyclic-flats-lattice}. Their proofs only use that $\mL$ is a finite modular lattice and that $\rho$ is monotone and submodular. Indeed, besides \Cref{lem:cover-transfer}, they only use parts (1), (2) and (4) of \Cref{lem:closure}, and these parts only need monotonicity and submodularity. For (1), if $\rho(X\vee K_1)=\rho(X\vee K_2)=\rho(X)$, then submodularity and monotonicity give $\rho(X\vee K_1\vee K_2)\leq\rho(X\vee K_1)+\rho(X\vee K_2)-\rho\bigl((X\vee K_1)\wedge(X\vee K_2)\bigr)\leq\rho(X)$.~Parts (2) and (4) follow from (1) and submodularity, as for matroids. So, with the same definitions of closure, flats and cyclic elements, these results hold for every real-valued monotone submodular function on a finite modular lattice. However, the assumption is essential for the reconstruction results with cyclic flats. Indeed, let $\mL=\{0_\mL<1_\mL\}$, and let $\rho(0_\mL)=0$ and $\rho(1_\mL)=\alpha$, where $0<\alpha<1$. Then $L=(\rho,\Ht,\mL)$ is an $\mathbb{R}$-latroid. Both elements of $\mL$ are flats, but the only cyclic flat of $L$ is $0_\mL$, of rank zero, independently of $\alpha$. Thus the family of flats alone does not determine~$\rho$. Even the cyclic flats together with their ranks do not determine~$\rho$. Moreover, we have
    $\nu(1_\mL)=1-\alpha>0=\nu(\cyc(1_\mL))$,
    so the identity $\nu(X)=\nu(\cyc(X))$ also fails.
\end{remark}

\begin{remark}\label{rem:chain-rings}
    Submodule lattices give further finite modular lattices that are neither distributive nor complemented, to which the results of this section apply, but neither \cite{byrne2025cyclic} nor \Cref{prop:convolution}. Let $R$ be a finite chain ring that is not a field, for instance $R=\ZZ/p^k\ZZ$ with $p$ prime and $k\geq2$, let $\mathfrak m$ be its maximal ideal, and let $n\geq2$. The lattice of submodules of $R^n$ is modular, as for every module, and its height function is the composition length. Let $S$ be the socle of $R^n$, that is, the sum of its simple submodules. Then $S$ is an $R/\mathfrak m$-vector space of dimension $n\geq2$, so the interval $[0,S]$ contains a copy of the diamond $M_3$, and the lattice is not distributive. Moreover, every nonzero submodule of $R^n$ contains a simple submodule, and hence meets $S$ nontrivially. Since $S\neq R^n$, the element $S$ has no complement. For $R=\ZZ/4\ZZ$, this is the subgroup lattice of $(\ZZ/4\ZZ)^n$, with $S=\{0,2\}^n$; the three subgroups of order $2$ of $\{0,2\}^2\times\{0\}^{n-2}$ give a copy of $M_3$. Lattices of this kind arise from codes over modules, which motivated the definition of latroids in \cite{vertigan2004latroids}.
\end{remark}

\section{$\ZZ$-Latroids on principal intervals of Young's lattice}\label{sec:young}

In this section we investigate $\ZZ$-latroids on principal intervals of Young's lattice. Throughout this section, we fix a partition
$$
    P_n=(\lambda_1,\ldots,\lambda_t)\vdash n,$$
with $\lambda_1\geq\cdots\geq\lambda_t>0$,
and consider the interval $\YL(P_n)=[\emptyset,P_n]\subseteq\YL$. This is a finite distributive, modular lattice. We set $\mL=\YL(P_n)$, $A=\ZZ$ and $\lVert P\rVert=|P|$. Since $|P|$ is the height of $P$ in $\YL(P_n)$, a function $\rho\colon\YL(P_n)\to\ZZ$ defines a $\ZZ$-latroid $L=(\rho,|\cdot|,\YL(P_n))$ if and only if $\rho(\emptyset)=0$ and $\rho$ satisfies $(L4)$ and $(L5)$. Moreover, the flats and the cyclic flats of $L$ form lattices, and $L$ is determined by its flats, and also by its cyclic flats together with their ranks, by the results in Section~\ref{sec:flats-cyclic}.

We write $\mI_L$, $\mB_L$ and $\mC_L$ for the collections of independent partitions, bases and circuits of~$L$.

Let $P\in\YL(P_n)$. A partition $P'\in\YL(P_n)$ covers $P$ if and only if $P'$ is obtained from $P$ by adding one box. We call a box $b\in P_n$ an \textbf{addable box} of $P$ if $b\notin P$ and $P\cup\{b\}$ is a Young diagram. Similarly, a box $b\in P$ is a \textbf{removable box} of $P$ if $P\setminus\{b\}$ is a Young diagram. For instance, let $P_n=(3,2)$ and $P=(2,1)$. In the picture below, the boxes of $P$ are shaded, its removable boxes are marked with $\mathsf{r}$, and its addable boxes are marked with $\mathsf{a}$:
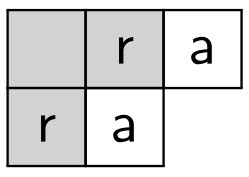
$$
    \ytableausetup{boxsize=1.3em}
    \begin{ytableau}
        *(gray!35) & *(gray!35)\mathsf{r} & \mathsf{a}\\
        *(gray!35)\mathsf{r} & \mathsf{a}
    \end{ytableau}
    \ytableausetup{smalltableaux}
$$
So $P$ is covered by $(3,1)$ and $(2,2)$, and it covers $(2)$ and $(1,1)$. By \Cref{lem:closure}\,(3) and \Cref{def:cyclic}, we can describe flats and cyclic partitions with these boxes:
\begin{enumerate}
    \item[(1)] $P$ is a flat if and only if adding each addable box of $P$ increases $\rho$.
    \item[(2)] $P$ is cyclic if and only if removing each removable box of $P$ does not change $\rho$.
\end{enumerate}

 We start with the modular rank functions, which count the boxes of a partition that lie in a fixed set of boxes.~Two rank functions that come from classical statistics of partitions, namely the size of the \emph{maximal hook} and the \emph{Durfee rank} of a partition, are of this form, and we describe them in detail. We finally construct paving $\ZZ$-latroids.

\subsection{Modular rank functions}

In this part we consider rank functions that count the boxes of $P$ that lie in a fixed set $W$ of boxes of $P_n$. We show that they are exactly the modular rank functions on $\YL(P_n)$. We also show how their independent partitions, bases, circuits, flats and cyclic flats can be read from~$W$.

Note that the removable boxes of a partition $P$ are its maximal boxes with respect to $\preceq$.

\begin{proposition}\label{prop:modular-classification}
    For a set of boxes $W\subseteq P_n$, let $\rho_W(P):=|P\cap W|$ for $P\in\YL(P_n)$. Then $L_W:=(\rho_W,|\cdot|,\YL(P_n))$ is a $\ZZ$-latroid, and $\rho_W$ is modular. Conversely, if $L=(\rho,|\cdot|,\YL(P_n))$ is a $\ZZ$-latroid and $\rho$ is modular, then $\rho=\rho_W$ for a unique $W\subseteq P_n$.
\end{proposition}

\begin{proof}
    Clearly $\rho_W(\emptyset)=0$. If $P\leq Q$, then $\rho_W(Q)-\rho_W(P)=|(Q\setminus P)\cap W|$ lies between $0$ and $|Q\setminus P|=|Q|-|P|$, so $(L4)$ holds. Since the join and the meet of $P$ and $Q$ in $\YL(P_n)$ are the union and the intersection of their Young diagrams, the inclusion--exclusion principle gives $\rho_W(P\vee Q)+\rho_W(P\wedge Q)=\rho_W(P)+\rho_W(Q)$. So $\rho_W$ is modular, and in particular $(L5)$ holds.

    Conversely, let $\rho$ be modular. For a box $b\in P_n$, the box $b$ is removable from $R_b$, and we let $w(b):=\rho(R_b)-\rho(R_b\setminus\{b\})$. Let $P\in\YL(P_n)$, and let $b$ be an addable box of $P$. Since $P\cup\{b\}$ is a partition containing $b$, we have $R_b\subseteq P\cup\{b\}$, and hence $R_b\setminus\{b\}\subseteq P$. So $R_b\vee P=P\cup\{b\}$ and $R_b\wedge P=R_b\setminus\{b\}$, and modularity gives
    $$
        \rho(P\cup\{b\})-\rho(P)=\rho(R_b)-\rho(R_b\setminus\{b\})=w(b).
    $$
    So adding the box $b$ always changes $\rho$ by $w(b)$. Adding the boxes of $P$ one at a time, starting from $\emptyset$, we obtain $\rho(P)=\sum_{b\in P}w(b)$. By \eqref{eq:unit-increase}, $w(b)\in\{0,1\}$ for every $b$, so $\rho=\rho_W$ with $W:=\{b\in P_n\mid w(b)=1\}$. Finally, $W$ is determined by $\rho$, since $b\in W$ if and only if $\rho(R_b)\neq\rho(R_b\setminus\{b\})$.
\end{proof}

In particular, there are exactly $2^n$ modular $\ZZ$-latroids on $\YL(P_n)$, one for each set of boxes of~$P_n$. The converse part of \Cref{prop:modular-classification} is a special case of a classical fact on distributive lattices that can be found in \cite{geissinger1973valuations}. Namely, a normalised modular function on a finite distributive lattice, that is, one with value $0$ at the least element, is determined by its values on the join-irreducible elements, that is, the elements that cover exactly one element. In $\YL(P_n)$, these elements are the rectangles $R_b$.

The structure of $L_W$ can be read from $W$ as follows.

\begin{proposition}\label{prop:modular-structure}
    Let $W\subseteq P_n$. Then the $\ZZ$-latroid $L_W$ has the following properties.
    \begin{enumerate}
        \item[(a)] A partition $P\in\YL(P_n)$ is independent if and only if $P\subseteq W$.
        \item[(b)] The rank of $L_W$ is $|W|$. It has a basis if and only if $W$ is a Young diagram, and in this case $W$ is its unique basis.
        \item[(c)] The circuits of $L_W$ are the rectangles $R_b$, where $b$ is a minimal box of $P_n\setminus W$ with respect to $\preceq$.
        \item[(d)] A partition $F\in\YL(P_n)$ is a flat if and only if all addable boxes of $F$ lie in $W$.
        \item[(e)] A partition $C\in\YL(P_n)$ is cyclic if and only if no removable box of $C$ lies in $W$.
    \end{enumerate}
    In particular, the cyclic flats of $L_W$ are the partitions whose addable boxes all lie in $W$ and whose removable boxes all lie outside $W$.
\end{proposition}

\begin{proof}
    (a) We have $\rho_W(P)=|P|$ if and only if $P\cap W=P$.

    (b) The rank is $\rho_W(P_n)=|W|$. By (a), a basis is a partition $B\subseteq W$ with $|B|=|W|$, that is, $B=W$.

    (c) By (a), a partition $P$ is dependent if and only if it contains a box $b\in P_n\setminus W$, that is, if and only if $R_b\subseteq P$ for some $b\in P_n\setminus W$. Hence the circuits are the minimal elements of $\{R_b\mid b\in P_n\setminus W\}$. Since $R_b\subseteq R_{b'}$ if and only if $b\preceq b'$, these are the rectangles $R_b$ with $b$ minimal in $P_n\setminus W$.

    (d) and (e) Adding or removing a box $b$ changes $\rho_W$ by one if $b\in W$, and does not change it otherwise. So both claims follow from the description of flats and cyclic partitions.
\end{proof}

\subsection{Maximal hook length}

In this part we show that the size of the maximal hook of a partition is a modular rank function. Afterwards we discuss the corresponding independent partitions, bases, circuits, flats and cyclic flats.

\begin{proposition}\label{prop:max-hook-rank}
    Let $W$ be the set of boxes of $P_n$ that lie in the first row or in the first column. The function $\rho_H:\YL(P_n)\rightarrow\ZZ$ defined by
    $\rho_H(P):=|H^P_{(1,1)}|$
    for $P\neq\emptyset$, together with $\rho_H(\emptyset):=0$, is equal to $\rho_W$. Thus
    $L_H=(\rho_H,|\cdot|,\YL(P_n))$
    is a $\ZZ$-latroid on $\YL(P_n)$, and $\rho_H$ is modular.
\end{proposition}

\begin{proof}
    For a nonempty partition $P\in\YL(P_n)$, the maximal hook $H^P_{(1,1)}$ consists of the boxes of $P$ in the first row or in the first column, that is, $H^P_{(1,1)}=P\cap W$. Hence $\rho_H=\rho_W$, and the claim follows from \Cref{prop:modular-classification}.
\end{proof}

In particular, adding a box to a partition increases $\rho_H$ by one if the box lies in the first row or in the first column. Otherwise, $\rho_H$ does not change. The same holds for removing a box.

\begin{example}\label{ex:max-hook-32}
    Let $n=5$ and $P_n=(3,2)$. Then $\YL(P_n)$ is the lattice of \Cref{fig:young-32}. For the $\ZZ$-latroid $L_H=(\rho_H,|\cdot|,\YL(P_n))$, the rank is $\rho_H(P_n)=4$, and $\mI_{L_H}=\YL(P_n)\setminus\{(2,2),(3,2)\}$, $\mB_{L_H}=\{(3,1)\}$, $\mC_{L_H}=\{(2,2)\}$, $\mF_{L_H}=\YL(P_n)\setminus\{(2,1),(3,1)\}$ and $\mZ_{L_H}=\{\emptyset,(2,2)\}$.
    The two cyclic flats have ranks $0$ and $3$. By \Cref{thm:reconstruction}, they determine the rank function. For instance, $\rho_H(P_n)=\min\{0+5,\,3+5-4\}=4$.
\end{example}

The next results determine the collections of independent partitions, bases, circuits, flats and cyclic flats of $L_H$ in general. We call a partition of the form $(a,1^{s-1})$, with $a,s\geq1$, a \textbf{hook partition}. In other words, all its rows except possibly the first have length one. A nonempty partition $P$ is a hook partition if and only if $H^P_{(1,1)}=P$.

\begin{corollary}\label{cor:max-hook-structure}
    Let $P_n=(\lambda_1,\ldots,\lambda_t)$ be nonempty. Then the following hold.
    \begin{enumerate}
        \item[(a)] The independent partitions are
        $$
            \mI_{L_H}
            =
            \{\emptyset\}
            \cup
            \{P\in\YL(P_n)
            \mid P\text{ is a hook partition}\}.
        $$
        \item[(b)] $L_H$ has the unique basis
        $$
            B_H=(\lambda_1,1^{t-1})=H_{(1,1)}^{P_n}.
        $$
        \item[(c)] The circuits are
        $$
            \mC_{L_H}
            =
            \begin{cases}
                \{(2,2)\},
                & \text{if } (2,2)\leq P_n,\\
                \emptyset,
                & \text{otherwise}.
            \end{cases}
        $$
    \end{enumerate}
\end{corollary}

\begin{proof}
    We apply \Cref{prop:modular-structure} with $W=H^{P_n}_{(1,1)}$.

    (a) A nonempty partition is contained in $W$ if and only if all its rows except possibly the first have length one, that is, if and only if it is a hook partition.

    (b) The set $W=(\lambda_1,1^{t-1})$ is a Young diagram.

    (c) If the box $(2,2)$ belongs to $P_n$, then it is the only minimal box of $P_n\setminus W$, since every box outside the first row and the first column lies weakly below and weakly to the right of $(2,2)$. So the only circuit is $R_{(2,2)}$, which is the partition $(2,2)$. Otherwise, $P_n$ is a hook partition, so $P_n\setminus W=\emptyset$ and there are no circuits.
\end{proof}

Finally we turn to the flats and the cyclic flats of $L_H$.

\begin{proposition}\label{prop:max-hook-flats}
    The empty partition is a flat of $L_H$. Let $F=(\mu_1,\ldots,\mu_s)\in\YL(P_n)$ be nonempty. Then $F$ is a flat of $L_H$ if and only if, for every $i\in\{2,\ldots,s\}$ such that $\mu_i<\mu_{i-1}$, one has $\mu_i=\lambda_i$.
\end{proposition}

\begin{proof}
    If $P_n=\emptyset$, then $\emptyset$ is the top element of $\YL(P_n)$, and hence a flat. Otherwise, the empty partition has only one cover, namely $(1)$, and $\rho_H((1))=1>\rho_H(\emptyset)=0$. Hence it is a flat. Now let $F=(\mu_1,\ldots,\mu_s)$ be nonempty. As noted after \Cref{prop:max-hook-rank}, adding a box to $F$ increases $\rho_H$ if the box lies in the first row or in the first column. This happens when the box extends the first row or creates a new row. If the added box lies in an existing row $i\geq2$, then $\rho_H$ does not change.
    A box can be added to an existing row $i\geq2$ precisely when $\mu_i<\mu_{i-1}$ and $\mu_i<\lambda_i$. Indeed, the first inequality guarantees that adding a box to row $i$ still gives a partition, whereas the second guarantees that the resulting partition is still contained in $P_n$. Consequently, $F$ is a flat if and only if no such box can be added. Equivalently, for every $i\in\{2,\ldots,s\}$ such that $\mu_i<\mu_{i-1}$, one must have $\mu_i=\lambda_i$.
\end{proof}

\begin{corollary}\label{cor:max-hook-cyclic}
    A nonempty partition $C=(\mu_1,\ldots,\mu_s)\in\YL(P_n)$ is cyclic in $L_H$ if and only if $s\geq2$, $\mu_1=\mu_2$ and $\mu_s\geq2$. Consequently, the cyclic flats of $L_H$ are $\emptyset$ and the partitions $C=(\mu_1,\ldots,\mu_s)\in\YL(P_n)$ such that
    \begin{enumerate}
        \item[(1)] $s\geq2$, $\mu_1=\mu_2$ and $\mu_s\geq2$, and
        \item[(2)] $\mu_i=\lambda_i$ for every $i\in\{3,\ldots,s\}$ with $\mu_i<\mu_{i-1}$.
    \end{enumerate}
\end{corollary}

\begin{proof}
    As noted after \Cref{prop:max-hook-rank}, removing a box from $C$ does not change $\rho_H$ if and only if the box lies neither in the first row nor in the first column. So $C$ is cyclic if and only if none of its removable boxes lies in the first row or in the first column. The removable boxes of $C$ are the last boxes $(i,\mu_i)$ of the rows $i$ with $\mu_i>\mu_{i+1}$, where we let $\mu_{s+1}:=0$. We look at them one by one.
    \begin{itemize}
        \item The box $(1,\mu_1)$ lies in the first row. It is removable if and only if $s=1$ or $\mu_1>\mu_2$.
        \item The box $(s,\mu_s)$ is always removable. It lies in the first row or in the first column if and only if $s=1$ or $\mu_s=1$.
        \item Every other removable box $(i,\mu_i)$ has $2\leq i\leq s-1$ and $\mu_i>\mu_{i+1}\geq1$. So it lies neither in the first row nor in the first column.
    \end{itemize}
    This proves the first claim. The second claim follows from \Cref{prop:max-hook-flats}, since $\emptyset$ is a cyclic flat and, when $\mu_1=\mu_2$, the condition of \Cref{prop:max-hook-flats} is empty for $i=2$.
\end{proof}

\begin{example}
    Let $P_n$ be a hook partition, for instance $P_n=(n)$ or $P_n=(1^n)$. Then every nonempty $P\in\YL(P_n)$ is a hook partition, so $\rho_H(P)=|P|$. Hence $L_H=(\rho_H,|\cdot|,\YL(P_n))$ is equal to the uniform $\ZZ$-latroid $\mU_n(\YL(P_n),\ZZ)$, also called the \textbf{free} $\ZZ$-latroid on $\YL(P_n)$. Moreover, we have $\mI_{L_H}=\mF_{L_H}=\YL(P_n)$, $\mB_{L_H}=\{P_n\}$, $\mC_{L_H}=\emptyset$ and $\mZ_{L_H}=\{\emptyset\}$.
\end{example}

\subsection{Durfee rank} 
In this part we consider the side length of the largest square contained in a Young diagram, that is, of its Durfee square. We call it the \emph{Durfee rank} of the partition. It is again a modular rank function, and we determine the corresponding independent partitions, bases, circuits, flats and cyclic flats.

\begin{definition}
    Let $P\in\YL(P_n)$. The \textbf{Durfee rank} of $P$, denoted by $\rho_D(P)$, is the side length of the largest square contained in $P$. Equivalently, if $P=(\mu_1,\ldots,\mu_s)$ is nonempty, then
    $$
        \rho_D(P)=\max\{1\leq i\leq s\mid \mu_i\geq i\}.
    $$
     We set $\rho_D(\emptyset):=0$.
\end{definition}

Recall that the boxes $(i,i)$ of a partition are its diagonal boxes. A partition $P$ contains $(i,i)$ if and only if $\mu_i\geq i$. Moreover, if $i\geq2$ and $\mu_i\geq i$, then also $\mu_{i-1}\geq\mu_i>i-1$. So the diagonal boxes of $P$ are exactly $(1,1),(2,2),\ldots,(r,r)$ with $r=\rho_D(P)$. In other words, $\rho_D(P)$ is the number of diagonal boxes of $P$. For example, $\rho_D((3,2))=2$, since $(3,2)$ contains $(1,1)$ and $(2,2)$.

\begin{proposition}\label{prop:durfee-rank}
    Let $W$ be the set of diagonal boxes of $P_n$. Then $\rho_D=\rho_W$ on $\YL(P_n)$. Hence $L_D=(\rho_D,|\cdot|,\YL(P_n))$ is a $\ZZ$-latroid on $\YL(P_n)$, and $\rho_D$ is modular.
\end{proposition}

\begin{proof}
    For $P\in\YL(P_n)$, the number $\rho_D(P)$ is the number of diagonal boxes of $P$, and these are exactly the boxes in $P\cap W$. Hence $\rho_D=\rho_W$, and the claim follows from \Cref{prop:modular-classification}.
\end{proof}

In particular, adding or removing a box changes $\rho_D$ by one if the box is a diagonal box, and does not change $\rho_D$ otherwise.

\begin{example}\label{ex:durfee-32}
    Let $P_n=(3,2)$, so that $\YL(P_n)$ is the lattice of \Cref{fig:young-32}. The $\ZZ$-latroid $L_D$ has rank $2$, and $\mI_{L_D}=\{\emptyset,(1)\}$, $\mB_{L_D}=\emptyset$, $\mC_{L_D}=\{(2),(1,1)\}$, $\mF_{L_D}=\{\emptyset,(3,1),(3,2)\}$ and $\mZ_{L_D}=\{\emptyset,(3,1)\}$.
    Thus, with the notion of basis used for latroids, $L_D$ is an example of a $\ZZ$-latroid having no basis. Moreover, its lattice of flats is the chain $\emptyset<(3,1)<(3,2)$, so it is not geometric.
\end{example}

The next result determines the independent partitions, bases and circuits of $L_D$ in general.

\begin{corollary}\label{cor:durfee-structure}
    Let $P_n$ be nonempty, and let $d:=\rho_D(P_n)$. Then the following hold.
    \begin{enumerate}
        \item[(a)] The independent partitions are
        $$
            \mI_{L_D}
            =
            \{\emptyset,(1)\}.
        $$
        \item[(b)] The bases are
        $$
            \mB_{L_D}
            =
            \begin{cases}
                \{(1)\},
                &\text{if }d=1,\\
                \emptyset,
                &\text{if }d\geq2.
            \end{cases}
        $$
        \item[(c)] The circuits are
        $$
            \mC_{L_D}
            =
            \{(2),(1,1)\}
            \cap\YL(P_n).
        $$
    \end{enumerate}
\end{corollary}

\begin{proof}
    We apply \Cref{prop:modular-structure} with $W=\{(1,1),\ldots,(d,d)\}$. 
    
    (a) A partition with at least two boxes contains $(1,2)$ or $(2,1)$, so the only partitions contained in $W$ are $\emptyset$ and $(1)$. 
    
    (b) The set $W$ is a Young diagram if and only if $d=1$, and then $W=(1)$. 
    
    (c) Every box $(i,j)$ of $P_n\setminus W$ satisfies $i<j$ or $i>j$, and hence $(1,2)\preceq(i,j)$ or $(2,1)\preceq(i,j)$. So the minimal boxes of $P_n\setminus W$ are those of $(1,2)$ and $(2,1)$ that lie in $P_n$, and the circuits are $R_{(1,2)}=(2)$ and $R_{(2,1)}=(1,1)$, whenever they lie in $\YL(P_n)$.
\end{proof}

We now determine the flats and the cyclic flats of $L_D$.

\begin{proposition}\label{prop:durfee-flats}
    Let $d:=\rho_D(P_n)$. For $0\leq j\leq d$, define $F_0:=\emptyset$ and, for $j\geq1$,
    $$
        F_j
        :=
        \bigcup_{i=1}^{j} H_{(i,i)}^{P_n}.
    $$
    In words, $F_j$ consists of the boxes of $P_n$ that lie in the first $j$ rows or in the first $j$ columns. Then $\mF_{L_D}=\{F_0,F_1,\ldots,F_d\}$. More precisely, if $P\in\YL(P_n)$ has Durfee rank $j$, then $\cl(P)=F_j$. In particular, the flats of $L_D$ form a chain, and $F_d=P_n$.
\end{proposition}

\begin{proof}
    Let $P\in\YL(P_n)$ and let $j:=\rho_D(P)$. If $j=0$, then $P=\emptyset$, and $\cl(P)=\emptyset=F_0$, since every nonempty partition has Durfee rank at least $1$. Assume from now on that $j\geq1$. Recall that the diagonal boxes of a partition are $(1,1),\ldots,(r,r)$ for some $r$, and the diagonal boxes of $P\vee Q=P\cup Q$ are those of $P$ together with those of $Q$. So for every $Q\in\YL(P_n)$ we have $\rho_D(P\vee Q)=\max\{\rho_D(P),\rho_D(Q)\}$. Hence
    $\rho_D(P\vee Q)=\rho_D(P)$ if and only if $\rho_D(Q)\leq j$.
    By the definition of closure,
    $\cl(P)
        =
        \bigvee
        \{Q\in\YL(P_n)\mid \rho_D(Q)\leq j\}$.
    A partition of Durfee rank at most $j$ cannot contain a box $(a,b)$ with both $a\geq j+1$ and $b\geq j+1$, since then it would also contain the box $(j+1,j+1)$ and would have Durfee rank at least $j+1$. Thus every such partition is contained in the union of the first $j$ rows and the first $j$ columns of $P_n$, which is $F_j$. Conversely, $F_j$ itself has Durfee rank $j$, since it contains the square $(j^j)$ and not the box $(j+1,j+1)$. Therefore $\cl(P)=F_j$.
    Hence, a partition of Durfee rank $j$ is a flat if and only if it is equal to $F_j$. Finally, since $d=\rho_D(P_n)$, there is no box of $P_n$ lying simultaneously below row $d$ and to the right of column $d$. Hence $F_d=P_n$.
\end{proof}

\begin{corollary}\label{cor:durfee-cyclic}
    Let $d$ and $F_0,\ldots,F_d$ be as in \Cref{prop:durfee-flats}. For $1\leq j\leq d$, the flat $F_j$ is cyclic if and only if $(j+1,j)\in P_n$ or $(j,j+1)\in P_n$. In particular, the cyclic flats of $L_D$ form a chain containing $\emptyset$, and $L_D$ is Schubert.
\end{corollary}

\begin{proof}
    Removing a box changes $\rho_D$ only if the box is a diagonal box. So $F_j$ is cyclic if and only if none of its diagonal boxes is removable. For $i<j$, the box $(i,i)$ is not removable from $F_j$, since the box $(i,i+1)$ to its right also lies in $F_j$. The box $(j,j)$ is removable from $F_j$ if and only if neither $(j+1,j)$ nor $(j,j+1)$ lies in $F_j$. Since $F_j$ contains all boxes of $P_n$ in the first $j$ rows and in the first $j$ columns, this holds if and only if neither $(j+1,j)$ nor $(j,j+1)$ lies in $P_n$. The last claim follows from \Cref{prop:durfee-flats}, since $\emptyset$ is a cyclic flat.
\end{proof}

\subsection{Paving $\ZZ$-latroids}

We next introduce paving latroids and give a construction on principal intervals of Young's lattice. We show that this construction gives exactly the paving $\ZZ$-latroids of rank $k$ in which every circuit of size $k$ is a flat. These play the role of sparse paving matroids.

\begin{definition}
    Let $L=(\rho,\lVert\cdot\rVert,\mL)$ be an $A$-latroid on $\mL$. We call $L$ \textbf{paving} if for every circuit $C$ of $L$ it holds that $\lVert C\rVert\geq\rho(L)$.
\end{definition}

Now let us return to our setting $\mL=\YL(P_n)$, $A=\ZZ$ and $\lVert P\rVert=|P|$. We write
$$
  \YL(P_n)_k:=\{P\in\YL(P_n)\mid |P|=k\},
  \qquad \textnormal{ and } \qquad
  \YL(P_n)_{\leq k}:=\{P\in\YL(P_n)\mid |P|\leq k\}.
$$

\begin{proposition}\label{prop:paving-construction}
    Let $2\leq k\leq n$ be an integer and let $\mS\subseteq\YL(P_n)_k$ be such that $P_n\notin\mS$ and $|P\wedge P'|\leq k-2$ for all distinct $P,P'\in\mS$. Then the function $\rho_\mS:\YL(P_n)\rightarrow\ZZ$ defined by
    $$
        \rho_\mS(P)=\begin{cases}
            k-1\quad&\text{if }P\in\mS,\\
            \min\{k,|P|\}\quad&\text{otherwise}
        \end{cases}
    $$
    satisfies axioms $(L4)$ and $(L5)$. Thus $L_\mS=(\rho_\mS,|\cdot|,\YL(P_n))$ is a paving $\ZZ$-latroid on $\YL(P_n)$ of rank $k$.
\end{proposition}

\begin{proof}
    Clearly $\rho_\mS(\emptyset)=0$, and since $P_n\notin\mS$ and $n\geq k$, we have $\rho_\mS(P_n)=k$. So the rank of $L_\mS$ is $k$.

    \textbf{Axiom $(L4)$.} We first consider a covering relation $P\lessdot P'$. Then
$|P'|=|P|+1$. If $P'\in\mS$, then $|P|=k-1$ and $\rho_\mS(P')=k-1=\rho_\mS(P)$. If
$P\in\mS$, then $|P'|=k+1$ and $\rho_\mS(P')=k=\rho_\mS(P)+1$. In all other cases,
$\rho_\mS$ agrees with the function $P\mapsto\min\{k,|P|\}$ on $P$ and $P'$, and the
assertion is clear. Hence
$$
    0\leq\rho_\mS(P')-\rho_\mS(P)\leq1 \quad\text{for all } P\lessdot P'.
$$
Now let $P<Q$ in $\YL(P_n)$, and choose a maximal chain
$P=P_0\lessdot P_1\lessdot\cdots\lessdot P_m=Q$. Each cover adds exactly one box,
so $m=|Q|-|P|$. Applying the previous inequality to each cover $P_{i-1}\lessdot P_i$
gives
$$
    0\leq\rho_\mS(Q)-\rho_\mS(P)=\sum_{i=1}^{m}\bigl(\rho_\mS(P_i)-\rho_\mS(P_{i-1})\bigr)\leq m=|Q|-|P|,
$$
which is $(L4)$.

    \textbf{Axiom $(L5)$.} Let $u(P):=\min\{k,|P|\}$, so that $\rho_\mS$ is obtained from $u$ by subtracting $1$ at the elements of $\mS$. First note that $u$ is submodular. Indeed, we have $|P|+|P'|=|P\vee P'|+|P\wedge P'|$ and $|P\wedge P'|\leq|P|,|P'|\leq|P\vee P'|$, and the function $x\mapsto\min\{k,x\}$ is concave.

    Now let $P,P'\in\YL(P_n)$. If $P$ and $P'$ are comparable, then $\{P\vee P',P\wedge P'\}=\{P,P'\}$, and $(L5)$ holds with equality. So assume that they are incomparable. We distinguish three cases.
    \begin{itemize}
        \item Neither $P$ nor $P'$ lies in $\mS$. Then $(L5)$ follows from the submodularity of $u$, since subtracting a possible unit at $P\vee P'$ or $P\wedge P'$ can only make the right-hand side smaller.
        \item Both $P$ and $P'$ lie in $\mS$. Then $|P|=|P'|=k$ and, by hypothesis, $|P\wedge P'|\leq k-2$. Since $|P\vee P'|=2k-|P\wedge P'|>k$, neither $P\vee P'$ nor $P\wedge P'$ lies in $\mS$, and we obtain
        $$
            \rho_\mS(P)+\rho_\mS(P')=2k-2\geq k+|P\wedge P'|=\rho_\mS(P\vee P')+\rho_\mS(P\wedge P').
        $$
        \item Exactly one of them lies in $\mS$, say $P\in\mS$ and $P'\notin\mS$. Since $P$ and $P'$ are incomparable, we have $|P\wedge P'|<k<|P\vee P'|$ and $|P\wedge P'|<|P'|$. Hence neither $P\vee P'$ nor $P\wedge P'$ lies in $\mS$, and $\rho_\mS(P\vee P')+\rho_\mS(P\wedge P')=k+|P\wedge P'|$. If $|P'|\geq k$, then this is at most $2k-1=\rho_\mS(P)+\rho_\mS(P')$, since $|P\wedge P'|\leq k-1$. If $|P'|<k$, then it is at most $k-1+|P'|=\rho_\mS(P)+\rho_\mS(P')$, since $|P\wedge P'|\leq|P'|-1$.
    \end{itemize}
    Hence $(L5)$ holds.

    It remains to check that $L_\mS$ is paving. Every partition of size at most $k-1$ is independent, so every circuit has size at least $k$. Since the rank of $L_\mS$ is $k$, the latroid $L_\mS$ is paving.
\end{proof}

\begin{remark}\label{rem:k=2}
    In contrast to the matroid and $q$-matroid theory, in the case $k=2$ we have $|\mS|\leq 1$. Indeed, any two distinct partitions $P,P'\in\YL(P_n)_2$ satisfy $|P\wedge P'|\geq 1$, because both contain $(1)$, whereas the condition in \Cref{prop:paving-construction} would require $|P\wedge P'|\leq 0$. Thus the only possibilities are $\mS=\emptyset$ or $\mS=\{P\}$ for some $P\in\YL(P_n)_2$ with $P\neq P_n$.
\end{remark}

Since two distinct partitions $P,P'$ of size $k$ satisfy $|P\wedge P'|=k-1$ if and only if both cover $P\wedge P'$, the condition on $\mS$ says that no two elements of $\mS$ cover a common partition.

Now we characterise the associated collections of independent partitions, bases, circuits, flats and cyclic flats.

\begin{proposition}\label{prop:paving-structure}
    Let $2\leq k\leq n$, and let $\mS\subseteq\YL(P_n)_k$ be such that $P_n\notin\mS$ and $|P\wedge P'|\leq k-2$ for all distinct $P,P'\in\mS$. Then the $\ZZ$-latroid $L_\mS$ has the following properties.
    \begin{enumerate}
        \item[(a)] The independent partitions are the elements of $\YL(P_n)_{\leq k}\setminus\mS$.

        \item[(b)] The bases are the elements of $\YL(P_n)_{k}\setminus\mS$.

        \item[(c)] The circuits are the elements of $\mS$ and the partitions $P\in\YL(P_n)_{k+1}$ that do not cover any element of $\mS$.

        \item[(d)] The flats are the elements of $\mS$, the partition $P_n$, and the partitions $P\in\YL(P_n)_{\leq k-1}$ that are not covered by any element of $\mS$.

        \item[(e)] The cyclic flats are $\emptyset$, the elements of $\mS$, and $P_n$ if $P_n$ is cyclic. Moreover, $P_n$ is cyclic if and only if $n\geq k+2$, or $n=k+1$ and $P_n$ does not cover any element of $\mS$.
    \end{enumerate}
    In particular, $L_\mS$ is Schubert if and only if $|\mS|\leq1$.
\end{proposition}

\begin{proof}
    A partition $P$ is independent if and only if $\rho_\mS(P)=|P|$. If $|P|<k$, this always holds. If $|P|=k$, it holds precisely when $P\notin\mS$, while every partition of size larger than $k$ is dependent. This proves (a). Since the rank of $L_\mS$ is $k$, statement (b) follows immediately.

    Every $P\in\mS$ is a circuit, because it is dependent and all partitions strictly contained in it have size at most $k-1$ and are independent. Let now $P\in\YL(P_n)_{k+1}$. Then $P$ is dependent, and it is a circuit if and only if all partitions of size $k$ contained in $P$ are independent. By (a), this is equivalent to saying that $P$ covers no element of $\mS$. Finally, a partition of size at least $k+2$ contains a dependent partition of size $k+1$, and therefore it is not a circuit. This proves (c).

    For (d), every element of $\mS$ is a flat, since each of its covers has rank $k$, whereas its rank is $k-1$. The top element $P_n$ is always a flat. If $|P|\leq k-2$, every cover of $P$ has rank $|P|+1$, and hence $P$ is a flat. If $|P|=k-1$, a cover $P'$ has the same rank as $P$ precisely when $P'\in\mS$. Hence such a $P$ is a flat if and only if it is not covered by an element of $\mS$. Finally, if $k\leq |P|<n$ and $P\notin\mS$, then $\rho_\mS(P)=k$ and every cover of $P$ also has rank $k$, so $P$ is not a flat.

    For (e), a nonempty partition $P\notin\mS$ with $|P|\leq k$ is not cyclic, since every partition covered by $P$ has rank $|P|-1=\rho_\mS(P)-1$. Every element of $\mS$ is cyclic, since the partitions it covers have rank $k-1$. A partition of size at least $k+2$ is cyclic, since all partitions it covers have rank $k$. A partition $P$ of size $k+1$ has rank $k$, and it is cyclic if and only if it covers no element of $\mS$. Moreover, $\emptyset$ is a flat by (d), since the elements of $\mS$ have size $k\geq2$. Together with (d), this proves (e). Since no element of $\mS$ contains another one, the last claim follows.
\end{proof}

The construction of \Cref{prop:paving-construction} gives exactly the following class of paving $\ZZ$-latroids.

\begin{theorem}\label{thm:paving-characterisation}
    Let $2\leq k\leq n$, and let $L=(\rho,|\cdot|,\YL(P_n))$ be a paving $\ZZ$-latroid of rank $k$. The following are equivalent.
    \begin{enumerate}
        \item[(a)] Every circuit $C$ of $L$ with $|C|=k$ is a flat of $L$.
        \item[(b)] $L=L_\mS$ for some $\mS\subseteq\YL(P_n)_k$ such that $P_n\notin\mS$ and $|P\wedge P'|\leq k-2$ for all distinct $P,P'\in\mS$.
    \end{enumerate}
    In this case, $\mS$ is the set of circuits of $L$ of size $k$.
\end{theorem}

\begin{proof}
    (b)$\Rightarrow$(a). By \Cref{prop:paving-structure}\,(c), the circuits of size $k$ of $L_\mS$ are exactly the elements of $\mS$, and they are flats by \Cref{prop:paving-structure}\,(d).

    (a)$\Rightarrow$(b). Let $\mS$ be the set of circuits of $L$ of size $k$. We show that $\rho=\rho_\mS$. Let $P\in\YL(P_n)$.

    If $|P|\leq k-1$, then $P$ is independent, that is, $\rho(P)=|P|$. Indeed, otherwise a minimal dependent partition contained in $P$ would be a circuit of size at most $k-1$, which is impossible since $L$ is paving.

    If $|P|=k$, then every partition strictly contained in $P$ is independent, so $P$ is dependent if and only if it is a circuit, that is, if and only if $P\in\mS$. Hence $\rho(P)=k$ if $P\notin\mS$. If $P\in\mS$, then $k-1=\rho(Q)\leq\rho(P)<k$ for any partition $Q$ covered by $P$, and so $\rho(P)=k-1$.

    If $|P|\geq k+1$, we show that $\rho(P)=k$ by induction on $|P|$. We have $\rho(P)\leq\rho(P_n)=k$. Let $Q$ be a partition covered by $P$. If $Q\in\mS$, then $Q$ is a flat by (a), and \Cref{lem:closure}\,(3) gives $\rho(P)>\rho(Q)=k-1$. Otherwise, $\rho(Q)=k$, by the previous case or by induction. In both cases $\rho(P)\geq k$.

    Thus $\rho=\rho_\mS$. It remains to check the assumptions on $\mS$. If $P_n\in\mS$, then $\rho(P_n)=k-1$, which contradicts the fact that $L$ has rank $k$. Let $P,P'\in\mS$ be distinct, and suppose that $|P\wedge P'|\geq k-1$. Then $|P\wedge P'|=k-1$ and $|P\vee P'|=k+1$. Since $\rho=\rho_\mS$, we get $\rho(P\wedge P')=k-1$ and $\rho(P\vee P')=k$, and submodularity gives the contradiction
    $$
        2k-2=\rho(P)+\rho(P')\geq\rho(P\vee P')+\rho(P\wedge P')=2k-1.
    $$
    Hence $|P\wedge P'|\leq k-2$.
\end{proof}

Following \cite{vertigan2004latroids,gorla2025latroids}, the \textbf{dual} of a $\ZZ$-latroid $L=(\rho,|\cdot|,\YL(P_n))$ is $L^*=(\rho^*,\lVert\cdot\rVert^*,\YL(P_n)^{\mathrm{op}})$, where
$\lVert X\rVert^*:=n-|X|$ and $\rho^*(X):=n-|X|-\rho(P_n)+\rho(X)$.

\begin{proposition}\label{prop:sparse-paving-dual}
    Let $2\leq k\leq n$, and let $L=(\rho,|\cdot|,\YL(P_n))$ be a paving $\ZZ$-latroid of rank $k$. The following are equivalent.
    \begin{enumerate}
        \item[(1)] Every circuit $C$ of $L$ with $|C|=k$ is a flat of $L$.
        \item[(2)] $\rho(X)=k$ for every $X\in\YL(P_n)$ with $|X|\geq k+1$.
        \item[(3)] The dual $L^*$ is paving.
    \end{enumerate}
\end{proposition}

\begin{proof}
    (1)$\Rightarrow$(2). By \Cref{thm:paving-characterisation}, $L=L_\mS$ for some $\mS\subseteq\YL(P_n)_k$, and $\rho_\mS(X)=k$ whenever $|X|\geq k+1$.

    (2)$\Rightarrow$(1). Let $C$ be a circuit with $|C|=k$. Every partition strictly contained in $C$ is independent, so $\rho(C)\geq k-1$, and $\rho(C)<k$ since $C$ is dependent. Hence $\rho(C)=k-1$. Every partition covering $C$ has size $k+1$, and hence rank $k$ by (2). So $C$ is a flat by \Cref{lem:closure}\,(3).

    (2)$\Leftrightarrow$(3). The greatest element of $\YL(P_n)^{\mathrm{op}}$ is $\emptyset$, so the rank of $L^*$ is $\rho^*(\emptyset)=n-k$. Moreover, $X$ is independent in $L^*$ if and only if $\rho^*(X)=\lVert X\rVert^*$, that is, if and only if $\rho(X)=k$.

    Suppose that (2) holds. Then every $X$ with $\lVert X\rVert^*\leq n-k-1$, that is, with $|X|\geq k+1$, is independent in $L^*$. So every circuit of $L^*$ has length at least $n-k$, and $L^*$ is paving.

    Conversely, suppose that $\rho(X)<k$ for some $X$ with $|X|\geq k+1$. Choose $Y\geq X$ in $\YL(P_n)$ maximal with $\rho(Y)<k$. Then $Y$ is dependent in $L^*$. The elements strictly below $Y$ in the reversed order are the partitions strictly containing $Y$. These have rank $k$ by the maximality of $Y$, so they are independent in $L^*$. Hence $Y$ is a circuit of $L^*$ with $\lVert Y\rVert^*=n-|Y|\leq n-|X|\leq n-k-1$, and $L^*$ is not paving.
\end{proof}

\begin{remark}\label{rem:sparse-paving}
    Recall that a \textbf{hyperplane} of a matroid of rank $r$ is a flat of rank $r-1$. A paving matroid of rank $r$ is called \textbf{sparse paving} if its dual matroid is also paving, or equivalently, if every circuit of size $r$ is a hyperplane. \Cref{prop:sparse-paving-dual} shows that both descriptions have a latroid version and that the two still agree. In a paving $\ZZ$-latroid of rank $k$, a circuit of size $k$ has rank $k-1$, so condition~(1) says that every such circuit is a flat of rank $k-1$, while condition~(3) says that the dual is paving. Here $L^*$ lives on the reversed lattice $\YL(P_n)^{\mathrm{op}}$, which is in general not isomorphic to a principal interval $\YL(Q)$ of Young's lattice. Indeed, for $Q\neq\emptyset$ the lattice $\YL(Q)$ has exactly one atom, whereas $\YL(P_n)^{\mathrm{op}}$ has one atom for each removable box of $P_n$. We therefore call a paving $\ZZ$-latroid of rank $k\geq 2$ \textbf{sparse paving} if it satisfies the equivalent conditions of \Cref{prop:sparse-paving-dual}. With this terminology, \Cref{thm:paving-characterisation} says that \Cref{prop:paving-construction} produces exactly the sparse paving $\ZZ$-latroids of rank $k$, for $2\leq k\leq n$.

    The class of paving $\ZZ$-latroids is strictly larger. For example, for $P_n=(4)$, the rank function with values $0,1,1,1,2$ on $\emptyset<(1)<(2)<(3)<(4)$ defines a paving $\ZZ$-latroid of rank $2$ whose only circuit is $(2)$. Since $\rho((3))=\rho((2))$, the circuit $(2)$ is not a flat, so this $\ZZ$-latroid is not sparse paving.
\end{remark}

By \Cref{prop:sparse-paving-dual}, for $2\leq k\leq n$, the sparse paving $\ZZ$-latroids of rank $k$ are exactly the paving $\ZZ$-latroids of rank $k$ whose dual is paving, as for matroids.

\section{Ferrers diagram rank-metric codes and their latroids}\label{sec:ferrers}

Ferrers diagram rank-metric codes were introduced by Etzion and Silberstein \cite{etzion2009error} for the construction of subspace codes. In this section, given a partition $D$, we associate a $\ZZ$-latroid on $\YL(D)$ to every Ferrers diagram rank-metric code with diagram $D$.

\subsection{The latroid of a Ferrers diagram rank-metric code}

Throughout this section, unless stated otherwise, we fix a nonempty partition $D:=P_n=(\lambda_1,\ldots,\lambda_t)$, and we let $h:=\lambda_1$. We work over $\F_q$, but everything in this section holds over any field. Ferrers diagrams are aligned at the top left, so they are Young diagrams in our sense. For $P\in\YL(D)$, the \textbf{full Ferrers support space} of $P$ is the space
$$
    \mC_P:=\{M\in\F_q^{t\times h}\mid M_{i,j}=0\text{ if }(i,j)\notin P\}
$$
of all matrices supported on $P$. Its dimension is $|P|$. The \textbf{support} $\supp(M)$ of a matrix $M$ is the set of boxes $(i,j)$ with $M_{i,j}\neq0$, and the support of a set of matrices is the union of their supports. Recall that $\down S$ is the smallest Young diagram containing a finite set $S$ of boxes. Since the removable boxes of $P$ are its maximal boxes, for $S\subseteq P$ we have $\down S=P$ if and only if $S$ contains all removable boxes of $P$.

\begin{definition}\label{def:ferrers-code}
    A \textbf{Ferrers diagram rank-metric code} with diagram $D$ is an $\F_q$-linear subspace $C\leq\mC_D$. Its \textbf{minimum rank distance} is $d(C):=\min\{\rank(M)\mid M\in C,\ M\neq0\}$ if $C\neq0$. For $P\in\YL(D)$, the \textbf{shortening} of $C$ to $P$ is the subspace $C(P):=C\cap\mC_P$ of all codewords that vanish outside $P$. We let $\rho_C(P):=|P|-\dim(C(P))$ and $L_C:=(\rho_C,|\cdot|,\YL(D))$.
\end{definition}

Recall that the \textbf{vector matroid} of a matrix $G$ over a field has the set of column indices of $G$ as ground set, and the rank of a set of indices is the dimension of the span of the corresponding columns. Distinct indices are distinct elements of the ground set, even if the corresponding columns coincide. The rank function $r$ satisfies $r(\emptyset)=0$, $r(S)\leq r(S\cup\{b\})\leq r(S)+1$ and $r(S\cup T)+r(S\cap T)\leq r(S)+r(T)$ for all sets of indices $S,T$ and every index $b$; see e.g. \cite{oxley2011matroid}. Recall also that the \textbf{uniform matroid} $U_{k,m}$ on an $m$-element ground set has rank function $S\mapsto\min\{k,|S|\}$. We use the standard bilinear form $\langle M,N\rangle=\sum_{(i,j)\in D}M_{i,j}N_{i,j}$ on $\mC_D$, and we write $C^\perp$ for the dual of $C$ with respect to it. For $P\in\YL(D)$, let $\pi_P\colon\mC_D\to\mC_P$ be the projection that sets all entries outside $P$ to zero.

\begin{proposition}\label{prop:code-latroid}
    Let $C\leq\mC_D$ be a Ferrers diagram rank-metric code. Then $\rho_C(P)=\dim(\pi_P(C^\perp))$ for every $P\in\YL(D)$, and $L_C$ is a $\ZZ$-latroid of rank $|D|-\dim(C)$.
\end{proposition}

\begin{proof}
    Let $P\in\YL(D)$. A matrix $X\in\mC_P$ is orthogonal to $\pi_P(C^\perp)$ if and only if it is orthogonal to $C^\perp$, since $X$ vanishes outside $P$. So the orthogonal complement of $\pi_P(C^\perp)$ in $\mC_P$ is $C\cap\mC_P=C(P)$, and hence $\dim(\pi_P(C^\perp))=|P|-\dim(C(P))=\rho_C(P)$.
    Now let $G$ be a matrix whose rows form a basis of $C^\perp$, with columns indexed by the boxes of $D$. Then $\rho_C(P)$ is the rank of the columns of $G$ indexed by the boxes of $P$, that is, the rank function $r$ of the vector matroid of $G$, evaluated at $P$. Since the join and the meet of two partitions are the union and intersection of their Young diagrams, the properties of $r$ recalled above show that $\rho_C$ satisfies $\rho_C(\emptyset)=0$, $(L4)$ and $(L5)$. Finally, $\rho_C(D)=|D|-\dim(C)$.
\end{proof}

In other words, $L_C$ is the restriction to $\YL(D)$ of the vector matroid of~$C^\perp$ on the boxes of~$D$. This construction is an instance of the code latroid of \cite[Example~5.8]{vertigan2004latroids}, for the lattice of the spaces $\mC_P$, $P\in\YL(D)$, with the dimension as length function. The vector matroid of $C^\perp$ does not need to coincide with the matroid of \Cref{prop:convolution}, although their rank functions agree on the Young subdiagrams. For instance, let $D=(2)$ and $C=\langle(0,1)\rangle$. Then $\rho_C$ takes the values $0,1,1$ on $\emptyset<(1)<(2)$, and the matroid of \Cref{prop:convolution} is the uniform matroid $U_{1,2}$, whose cyclic flats are $\emptyset$ and $D$. The vector matroid of $C^\perp=\langle(1,0)\rangle$ has the coloop $(1,1)$ and the loop $(1,2)$, and its only cyclic flat is $\{(1,2)\}$.

We will often use the following observation. If $P'$ is obtained from $P$ by adding a box $b$, then $C(P)$ is the kernel of the linear map $C(P')\to\F_q$, $M\mapsto M_b$. Hence,
\begin{equation}\label{eq:code-cover}
    \begin{gathered}
        \dim(C(P'))-\dim(C(P))\in\{0,1\},\\
        \rho_C(P')-\rho_C(P)=1-\bigl(\dim(C(P'))-\dim(C(P))\bigr),
    \end{gathered}
\end{equation}
and $\dim(C(P'))=\dim(C(P))+1$ if and only if some codeword of $C(P')$ is nonzero at $b$.

\subsection{Cyclic flats and shortenings}

We now describe the structure of $L_C$ in terms of the shortenings of $C$.

\begin{proposition}\label{prop:code-dictionary}
    Let $C\leq\mC_D$ be a Ferrers diagram rank-metric code and let $P\in\YL(D)$. Then the following hold.
    \begin{enumerate}
        \item[(1)] $P$ is independent in $L_C$ if and only if $C(P)=0$.
        \item[(2)] $P$ is a circuit of $L_C$ if and only if $C(P)=\langle M\rangle$ for some nonzero matrix $M$ which is nonzero at every removable box of $P$.
        \item[(3)] $P$ is a flat of $L_C$ if and only if, for every addable box $a$ of $P$, every codeword supported on $P\cup\{a\}$ vanishes at $a$.
        \item[(4)] $P$ is cyclic if and only if every removable box of $P$ lies in $\supp(C(P))$. Moreover, $\cyc(P)=\down\supp(C(P))$.
    \end{enumerate}
\end{proposition}

\begin{proof}
    (1) By definition, $P$ is independent if and only if $\rho_C(P)=|P|$, that is, $C(P)=0$.

    (2) Suppose that $P$ is a circuit. Then $C(P)\neq0$, and $C(P\setminus\{b\})=0$ for every removable box $b$ of $P$. By \eqref{eq:code-cover}, $C(P)$ is one-dimensional, and a spanning matrix is nonzero at every removable box of $P$. Conversely, suppose that $C(P)=\langle M\rangle$ for such a matrix $M\neq0$. Then $C(P)\neq0$, so $P$ is dependent. Moreover, $M\notin C(P\setminus\{b\})$, and hence $C(P\setminus\{b\})=0$, for every removable box $b$ of $P$. Every partition strictly below $P$ is contained in some $P\setminus\{b\}$, so it is independent by (1). Hence $P$ is a circuit.

    (3) By \Cref{lem:closure}\,(3), $P$ is a flat if and only if $\rho_C(P\cup\{a\})=\rho_C(P)+1$ for every addable box $a$ of $P$. By \eqref{eq:code-cover}, this means that $C(P\cup\{a\})=C(P)$, that is, every codeword supported on $P\cup\{a\}$ vanishes at $a$.

    (4) By \Cref{def:cyclic} and \eqref{eq:code-cover}, $P$ is cyclic if and only if $\dim(C(P\setminus\{b\}))=\dim(C(P))-1$ for every removable box $b$ of $P$, that is, if and only if every removable box of $P$ lies in $\supp(C(P))$. Let $K:=\down\supp(C(P))$. By the observation before \Cref{def:ferrers-code}, $P$ is cyclic if and only if $K=P$. Since every codeword of $C(P)$ vanishes outside $K$, we have $C(K)=C(P)$, and hence $K$ is cyclic. If $Q\leq P$ is cyclic, then $Q=\down\supp(C(Q))\subseteq K$, since $C(Q)\subseteq C(P)$. So $K$ is the largest cyclic element below $P$, and $K=\cyc(P)$ by \Cref{lem:cyclic}\,(2).
\end{proof}

In (2), a single codeword must be nonzero at every removable box of $P$, whereas in (4) different codewords may be nonzero at different removable boxes.

Applying \Cref{thm:reconstruction} to $L_C$, we obtain the following result.

\begin{corollary}\label{thm:shortening-profile}
    Let $C\leq\mC_D$ be a Ferrers diagram rank-metric code. Then, for every $P\in\YL(D)$,
    $$
        \dim(C(P))=\max_{Z\in\mZ_{L_C}}\Bigl(\dim(C(Z))-|Z\setminus P|\Bigr),
    $$
    and the maximum is attained at $Z=\down\supp(C(\cl(P)))$. In particular, the cyclic flats of $L_C$, together with the dimensions of the corresponding shortenings, determine the dimensions of all shortenings of $C$.
\end{corollary}

\begin{proof}
    Since $|P\vee Z|=|P\cup Z|=|P|+|Z\setminus P|$, \Cref{thm:reconstruction} gives
    $$
        \dim(C(P))=|P|-\rho_C(P)=\max_{Z\in\mZ_{L_C}}\bigl(|P|-\rho_C(Z)-|P\cup Z|+|Z|\bigr)=\max_{Z\in\mZ_{L_C}}\bigl(\dim(C(Z))-|Z\setminus P|\bigr),
    $$
    and the maximum is attained at $\cyc(\cl(P))=\down\supp(C(\cl(P)))$, by \Cref{prop:code-dictionary}\,(4).
\end{proof}

\begin{example}\label{ex:ferrers-one-codeword}
    Let $D=(3,2)$, and let $C\leq\mC_D$ be the code spanned by
    $$
        M=\begin{pmatrix}1&0&0\\0&1&0\end{pmatrix}.
    $$
    Then $d(C)=2$. Moreover, $C(P)=C$ if $(2,2)\leq P$, and $C(P)=0$ otherwise. So $L_C$ is the maximal hook latroid $L_H$ of \Cref{ex:max-hook-32}. By \Cref{prop:code-dictionary}, its cyclic flats are $\emptyset$ and $(2,2)$, and $(2,2)$ is its only circuit. \Cref{thm:shortening-profile} reads $\dim(C(P))=\max\{0,\,1-|(2,2)\setminus P|\}$.
\end{example}

\subsection{The minimum rank distance}

The minimum rank distance of $C$ is related to $L_C$ through the Etzion--Silberstein bound; see \cite{etzion2009error}. Let $d\geq1$ and $P\in\YL(D)$. For $0\leq j\leq d-1$, let $\nu_j(P,d)$ be the number of boxes of $P$ that are neither in the first $j$ columns nor in the first $d-1-j$ rows, and let $\nu_{\min}(P,d):=\min_{0\leq j\leq d-1}\nu_j(P,d)$. These numbers should not be confused with the nullity $\nu$ of \Cref{sec:flats-cyclic}. Finally, let $T_d:=(d,d-1,\ldots,1)$.

\begin{theorem}[{\cite[Theorem~1]{etzion2009error}}]\label{thm:ES-bound}
    Let $C$ be a nonzero Ferrers diagram rank-metric code with diagram $P$ and minimum rank distance $d$. Then $\dim(C)\leq\nu_{\min}(P,d)$.
\end{theorem}

\begin{proposition}\label{prop:code-distance}
    Let $C\leq\mC_D$ be a nonzero Ferrers diagram rank-metric code with minimum rank distance $d$. Then the following hold.
    \begin{enumerate}
        \item[(1)] $\dim(C(P))\leq\nu_{\min}(P,d)$ for every $P\in\YL(D)$.
        \item[(2)] If $P\in\YL(D)$ satisfies $T_d\not\leq P$, then $P$ is independent in $L_C$. In particular, every circuit and every nonempty cyclic element of $L_C$ contains $T_d$.
    \end{enumerate}
\end{proposition}

\begin{proof}
    (1) If $C(P)=0$, the inequality holds since $\nu_{\min}(P,d)\geq0$. Otherwise, $C(P)$ is a Ferrers diagram rank-metric code with diagram $P$ and minimum rank distance $d'\geq d$, so $\dim(C(P))\leq\nu_{\min}(P,d')$ by \Cref{thm:ES-bound}. Moreover, $\nu_j(P,d+1)\leq\nu_j(P,d)$ for every $j\leq d-1$, since one more row is deleted. Hence $\nu_{\min}(P,d')\leq\nu_{\min}(P,d)$.

(2) Suppose that $T_d\not\leq P$, and write
$P=(\mu_1,\ldots,\mu_s)$, with $\mu_i:=0$ for $i>s$.
Then $\mu_i\leq d-i$ for some $1\leq i\leq d$.
Deleting the first $i-1$ rows and the first $d-i$ columns
therefore removes every box of $P$. Hence
$\nu_{d-i}(P,d)=0$, so $\nu_{\min}(P,d)=0$ (see also \cite[Lemma~4.17]{beeloo2026irreducible}).
By (1), $C(P)=0$, so $P$ is independent.
Finally, if $Z$ is a circuit or a nonempty cyclic element,
then $C(Z)\neq0$ by \Cref{prop:code-dictionary},
and hence $T_d\leq Z$.
\end{proof}

\begin{remark}\label{rem:code-latroid-scope}
    The latroid $L_C$ only depends on the supports of the codewords of $C$ in the fixed Ferrers coordinates. For instance, let $D=(2,2)$. The codes spanned by the identity matrix and by the all-one matrix have the same latroid, since in both cases $C(P)\neq0$ only for $P=D$, but their minimum rank distances are $2$ and $1$. Moreover, let $E_{i,j}$ be the matrix with a single $1$ in position $(i,j)$. The codes spanned by $E_{1,1}$ and by $E_{2,1}$ are equivalent, by exchanging the two rows, but the partition $(1)$ is a circuit of the first latroid and independent in the second one. So $L_C$ is not an invariant of the rank-metric equivalence class of $C$.
\end{remark}

\begin{proposition}\label{prop:not-from-code}
    For $D=(8,7,\ldots,1)$ there is a $\ZZ$-latroid on $\YL(D)$ which is not of the form $L_C$ for any code $C$ over any field.
\end{proposition}

\begin{proof}
    Let $D=(8,7,\ldots,1)$, let $R$ be the set of the $8$ removable boxes of $D$, and let $r$ be the rank function of the V\'amos matroid on $R$, an eight-element matroid that is not representable over any field; see \cite{oxley2011matroid}. Since $r$ is a matroid rank function and intersecting with $R$ preserves unions and intersections, $\rho(P):=r(P\cap R)$ defines a $\ZZ$-latroid on $\YL(D)$. Suppose that $\rho=\rho_C$ for some code $C$. By \Cref{prop:code-latroid} and its proof, there is a representable matroid $M$ on the boxes of $D$ with rank function $r_M$ such that $\rho(P)=r_M(P)$ for every $P\in\YL(D)$. Let $\overline{D}:=D\setminus R$. Recall that the \textbf{contraction} $M/\overline{D}$ is the matroid on the remaining boxes with rank function
    $$
        r_{M/\overline{D}}(S)=r_M(\overline{D}\cup S)-r_M(\overline{D}),
    $$
    and that contractions of matroids representable over a field are representable over the same field; see \cite{oxley2011matroid}. For every $S\subseteq R$, the set $\overline{D}\cup S$ is a Young diagram, and hence $r(S)=\rho(\overline{D}\cup S)-\rho(\overline{D})=r_{M/\overline{D}}(S)$. So the V\'amos matroid is the contraction $M/\overline{D}$, which is representable, a contradiction.
\end{proof}

\begin{remark}\label{rem:small-codes}
Over a fixed field $\F_q$, the proof of \Cref{prop:not-from-code} gives much smaller examples. Let $D=(q+2,q+1,\ldots,1)$, let $R$ be the set of its $q+2$ removable boxes, and let $r$ be the rank function of the uniform matroid $U_{2,q+2}$ on $R$. The matroid $U_{2,q+2}$ is not representable over $\F_q$, since a representation would consist of $q+2$ pairwise linearly independent vectors spanning a two-dimensional vector space over $\F_q$, which has only $q+1$ one-dimensional subspaces. Hence the same argument shows that the $\ZZ$-latroid $(\rho,|\cdot|,\YL(D))$ with $\rho(P):=r(P\cap R)$ is not of the form $L_C$ for any code $C$ over $\F_q$. For $q=2$, the diagram $D=(4,3,2,1)$ has $10$ boxes. We do not know the smallest diagram $D$ that carries a $\ZZ$-latroid which does not come from a code over any field.
\end{remark}

\section{$q$-Polymatroids from cyclic flats}\label{sec:qpoly}

We fix a partition $P_n=(\lambda_1,\ldots,\lambda_t)$, which we now assume to be nonempty, and consider $\ZZ$-latroids $L=(\rho,|\cdot|,\YL(P_n))$. We look at every cyclic flat of $L$ as a Ferrers diagram. To each such diagram we associate the space of all matrices supported on it, and we study the $q$-polymatroids defined by these spaces. In \Cref{sec:ferrers} we attached a $\ZZ$-latroid to a code, while here we attach matrix spaces, and hence $q$-polymatroids, to the cyclic flats of a $\ZZ$-latroid. Note that they only depend on the shapes of the cyclic flats, and not on their ranks. Inclusions between cyclic flats give inclusions between these spaces, and hence a rank-difference inequality, which gives a $q$-quotient relation when the associated $q$-polymatroids are $q$-matroids. 

We now recall some definitions on $q$-polymatroids. Let $q$ be a prime power, let $h:=\lambda_1\geq1$, and let $E=\F_q^t$. A \textbf{$q$-polymatroid} on $E$ is a function $r$ that assigns a real number to every subspace of $E$, such that $0\leq r(V)\leq\dim(V)$, $r(V)\leq r(V')$ whenever $V\leq V'$, and $r(V+V')+r(V\cap V')\leq r(V)+r(V')$ for all $V,V'$. A \textbf{$q$-matroid} is a $q$-polymatroid with integer values. A \textbf{flat} of a $q$-matroid is a subspace $F$ such that $r(F+\langle x\rangle)>r(F)$ for every $x\in E\setminus F$. The \textbf{space of loops} is the largest subspace of rank $0$. For $0\leq k\leq t$, the \textbf{uniform $q$-matroid} $\mU_{k,t}(q)$ has rank function $V\mapsto\min\{k,\dim(V)\}$. We will use the following two notions.

\begin{definition}
    Let $\mN_1$ and $\mN_2$ be $q$-matroids on the same space $E$. We say that $\mN_1$ is a \textbf{$q$-quotient} of $\mN_2$ if every flat of $\mN_1$ is a flat of $\mN_2$, that is, if the identity map $\mN_2\to\mN_1$ is a strong map; see \cite[Definition~2.21]{CJ25qanadeltamats}.
\end{definition}

\begin{definition}
    Let $\mT=\{\mN_1,\ldots,\mN_m\}$ be a set of distinct $q$-matroids on the same space $E$. We call $\mT$ a \textbf{$q$-quotient collection} if, for all $i\neq j$, $\mN_i$ is a $q$-quotient of $\mN_j$ or $\mN_j$ is a $q$-quotient of $\mN_i$.
\end{definition}

For matroids, quotients, which correspond to strong identity maps between matroids on the same ground set, are well-studied (see e.g. \cite{kung1986strong}), and a chain of quotients is a flag matroid; see~\cite{borovik2003coxeter}. So a $q$-quotient collection can be seen as a $q$-analogue of a flag matroid. Quotients appear naturally in coding theory. For matroids defined by generator matrices, the matroid of a subcode is a quotient of the matroid of the containing code. The generator-matrix construction also gives the corresponding quotient relation for nested $\F_{q^m}$-linear rank-metric codes; see \cite[Proposition~5.2]{CJ25qanadeltamats}. Recently, Ceria, Johnsen and Jurrius used $q$-quotients to define $q$-analogues of generalized matroids and of $\Delta$-matroids \cite{CJ25qanadeltamats}, following the classical constructions of \cite{tardos1985generalized,bouchet1987greedy}.

For every Young diagram $P\subseteq(h^t)$, consider the full Ferrers support space of $P$
$$
    \mC_P:=\{M\in\F_q^{t\times h}\mid M_{i,j}=0\text{ if }(i,j)\notin P\}
$$
of all matrices supported on $P$, which has dimension $|P|$. Since $h=\lambda_1$, every $P\in\YL(P_n)$ is such a diagram. We mostly use $\mC_P$ for $P\in\YL(P_n)$, but we also need the rectangles $(h^s)$ with $0\leq s\leq t$, which do not need to lie in $\YL(P_n)$. Let $e_1,\ldots,e_t$ be the standard basis of~$E$, and write $\colsp(M)\leq E$ for the column space of a matrix $M$. For a subspace $V\leq E$, let $\mC_P(V):=\{M\in\mC_P\mid\colsp(M)\leq V\}$, and let $V^\perp$ be the orthogonal complement of $V$ with respect to the standard inner product. In \cite{gorla2020rank} it is shown that the space $\mC_P$ defines a $q$-polymatroid $\mM_{\mC_P}$ on $E$ with rank function
$$
    r_{\mC_{P}}(V)=\frac{\dim_{\F_q}(\mC_P)-\dim_{\F_q}(\mC_P(V^\perp))}{h}.
$$
The authors of \cite{gorla2020rank} assume $t\leq h$. If $t>h$, the same construction is obtained from the row $q$-polymatroid of the transposed code, which still has ground space $\F_q^t$ and denominator $h$. In any case, \Cref{lem:average} below shows directly that $r_{\mC_{P}}$ is a $q$-polymatroid for all $t$ and $h$. We use the convention $(h^0)=\emptyset$.

\begin{proposition}\label{prop:qmatroid-rectangles}
    Let $P=(\mu_1,\ldots,\mu_t)\in\YL(P_n)$, where we add zero parts if necessary. Then $\mM_{\mC_P}$ is a $q$-matroid if and only if $P=(h^s)$ for some $0\leq s\leq t$.
\end{proposition}

\begin{proof}
    Since $\mC_P$ consists of all matrices supported on $P$, we have $\dim_{\F_q}(\mC_P)=|P|=\mu_1+\cdots+\mu_t$. In particular, $r_{\mC_{P}}(E)=|P|/h$.
    For $1\leq i\leq t$, let $H_i:=\langle e_1,\ldots,e_{i-1},e_{i+1},\ldots,e_t\rangle$, so that $H_i^\perp=\langle e_i\rangle$. The matrices in $\mC_P(\langle e_i\rangle)$ are the matrices in $\mC_P$ that are zero outside the $i$-th row. Since the $i$-th row of $P$ has $\mu_i$ boxes, we get $\dim_{\F_q}(\mC_P(\langle e_i\rangle))=\mu_i$. Therefore $r_{\mC_{P}}(H_i)=(|P|-\mu_i)/h$.

    Suppose that $\mM_{\mC_P}$ is a $q$-matroid. Then $r_{\mC_{P}}(E)$ and $r_{\mC_{P}}(H_i)$ are integers for every $i$, and hence so is $\mu_i/h=r_{\mC_{P}}(E)-r_{\mC_{P}}(H_i)$. Since $0\leq\mu_i\leq h$, we get $\mu_i\in\{0,h\}$ for every $i$. Since $P$ is a partition, $P=(h^s)$ for some $0\leq s\leq t$.

    Conversely, suppose that $P=(h^s)$. Let $U_s:=\langle e_1,\ldots,e_s\rangle\leq E$, with $U_0:=\{0\}$. Then $\mC_P$ consists of the matrices whose columns lie in $U_s$. Hence $\dim_{\F_q}(\mC_P)=sh$ and $\dim_{\F_q}(\mC_P(V^\perp))=h\dim(U_s\cap V^\perp)$ for every $V\leq E$. So $r_{\mC_{P}}(V)=s-\dim(U_s\cap V^\perp)=\dim(V)-\dim(V\cap U_s^\perp)$, where the second equality follows from Grassmann's identity. Thus $r_{\mC_{P}}$ only takes integer values, and $\mM_{\mC_P}$ is a $q$-matroid.
\end{proof}

We can describe $\mM_{\mC_P}$ for every $P$ by using the column lengths of $P$. For $0\leq s\leq t$, let $U_s:=\langle e_1,\ldots,e_s\rangle\leq E$, with $U_0=\{0\}$ as before, and let $\mM_s$ be the $q$-matroid on $E$ with rank function $r_s(V)=\dim(V)-\dim(V\cap U_s^\perp)$. By Grassmann's identity, $r_s(V)=\dim(V+U_s^\perp)-(t-s)$. Since $V\mapsto\dim(V+U_s^\perp)$ is monotone and submodular, and $0\leq r_s(V)\leq\dim(V)$, the function $r_s$ is indeed the rank function of a $q$-matroid, for all $t$ and $h$. Its space of loops is $U_s^\perp$, and $r_s(V)=\dim(V)$ for every $V\leq U_s$. By the proof of \Cref{prop:qmatroid-rectangles}, we have $\mM_s=\mM_{\mC_{(h^s)}}$.

\begin{lemma}\label{lem:average}
    Let $P\in\YL(P_n)$, and let $\mu'_1\geq\cdots\geq\mu'_h$ be the column lengths of $P$, that is, $\mu'_j=|\{i\mid (i,j)\in P\}|$. Then $r_{\mC_{P}}(V)=\frac1h\sum_{j=1}^{h}r_{\mu'_j}(V)$ for every $V\leq E$.
\end{lemma}

\begin{proof}
    The $j$-th column of a matrix $M\in\mC_P$ lies in $U_{\mu'_j}$, and the columns can be chosen independently. Moreover, $\colsp(M)\leq V^\perp$ if and only if every column of $M$ lies in $V^\perp$. Hence $\dim(\mC_P(V^\perp))=\sum_j\dim(U_{\mu'_j}\cap V^\perp)$, while $\dim(\mC_P)=\sum_j\mu'_j$. As in the proof of \Cref{prop:qmatroid-rectangles}, we have $\mu'_j-\dim(U_{\mu'_j}\cap V^\perp)=\dim(V)-\dim(V\cap U_{\mu'_j}^\perp)=r_{\mu'_j}(V)$, and the claim follows.
\end{proof}

So the rank function of $\mM_{\mC_P}$ is the average of the rank functions of $\mM_{\mu'_1},\ldots,\mM_{\mu'_h}$. Since the axioms of a $q$-polymatroid are preserved under taking averages, this shows directly that $\mM_{\mC_P}$ is a $q$-polymatroid, for all $t$ and $h$. In particular, it has integer values when all columns of $P$ have the same length, that is, when $P=(h^s)$.

The next result shows that inclusions of diagrams give the rank-difference condition that characterises $q$-quotients; see also \cite[Proposition~2.22]{CJ25qanadeltamats}.

\begin{proposition}\label{prop:rank-difference}
    Let $P\leq Q$ in $\YL(P_n)$. Then
    $$
        r_{\mC_{Q}}(X)-r_{\mC_{Q}}(Y)\geq r_{\mC_{P}}(X)-r_{\mC_{P}}(Y)\quad\text{for all }Y\leq X\leq E.
    $$
    In particular, if $\mM_{\mC_P}$ and $\mM_{\mC_Q}$ are $q$-matroids, then $\mM_{\mC_P}$ is a $q$-quotient of $\mM_{\mC_Q}$.
\end{proposition}

\begin{proof}
    We follow the proof of \cite[Proposition~5.2]{CJ25qanadeltamats}, which treats nested $\F_{q^m}$-linear codes; the same argument works for the $\F_q$-linear spaces $\mC_P\leq\mC_Q$. Let $Y\leq X\leq E$. Then $X^\perp\leq Y^\perp$, so $\mC_P(X^\perp)\leq\mC_P(Y^\perp)$ and $\mC_Q(X^\perp)\leq\mC_Q(Y^\perp)$. Since $P\leq Q$, we have $\mC_P\leq\mC_Q$, and $\mC_P(V)=\mC_P\cap\mC_Q(V)$ for every $V\leq E$. In particular, $\mC_P(Y^\perp)\cap\mC_Q(X^\perp)=\mC_P(X^\perp)$, so the linear map
    $$
        \mC_P(Y^\perp)/\mC_P(X^\perp)\longrightarrow\mC_Q(Y^\perp)/\mC_Q(X^\perp)
    $$
    induced by the inclusion is injective. By the definition of $r_{\mC_P}$ and $r_{\mC_Q}$, we get
    $$
        \begin{aligned}
            r_{\mC_P}(X)-r_{\mC_P}(Y)&=\frac{\dim(\mC_P(Y^\perp))-\dim(\mC_P(X^\perp))}{h}\\
            &\leq\frac{\dim(\mC_Q(Y^\perp))-\dim(\mC_Q(X^\perp))}{h}=r_{\mC_Q}(X)-r_{\mC_Q}(Y).
        \end{aligned}
    $$
    By \cite[Proposition~2.22]{CJ25qanadeltamats}, for $q$-matroids the first claim means that the identity map $\mM_{\mC_Q}\to\mM_{\mC_P}$ is a strong map, that is, $\mM_{\mC_P}$ is a $q$-quotient of $\mM_{\mC_Q}$.
\end{proof}

The proof only uses that $\mC_P\leq\mC_Q$, so the same inequality holds for any two nested spaces of $t\times h$ matrices. In particular, for every $\ZZ$-latroid $L$ and all $Z\leq Z'$ in $\mZ_L$, the $q$-polymatroids $\mM_{\mC_Z}$ and $\mM_{\mC_{Z'}}$ satisfy the rank-difference condition of \Cref{prop:rank-difference}. For Ferrers diagrams, the converse also holds: if $P,Q\in\YL(P_n)$ satisfy $r_{\mC_{Q}}(X)-r_{\mC_{Q}}(Y)\geq r_{\mC_{P}}(X)-r_{\mC_{P}}(Y)$ for all $Y\leq X\leq E$, then $P\leq Q$. Indeed, for $1\leq k\leq t$, take $X=U_k$ and $Y=U_{k-1}$. Then \Cref{lem:average} gives $r_{\mC_P}(U_k)-r_{\mC_P}(U_{k-1})=\mu_k/h$, where $\mu_k$ is the length of the $k$-th row of $P$, and similarly for $Q$. So the $k$-th row of $P$ is at most as long as the $k$-th row of $Q$ for every $k$, that is, $P\leq Q$. Hence $P\leq Q$ if and only if $\mM_{\mC_P}$ and $\mM_{\mC_Q}$ satisfy the rank-difference condition. This suggests studying analogues of $q$-quotients, $q$-g-matroids and $q$-$\Delta$-matroids for $q$-polymatroids.

\begin{example}
    Let $P_n=(2,2)$ and $q=2$, and list the elements of $\YL(P_n)$ in the order $\emptyset,(1),(1,1),(2),(2,1),(2,2)$. In this order, the spaces $\mC_P\leq\F_2^{2\times2}$ consist of the matrices of the form
    $$
        \begin{pmatrix}0&0\\0&0\end{pmatrix},\quad
        \begin{pmatrix}a&0\\0&0\end{pmatrix},\quad
        \begin{pmatrix}a&0\\b&0\end{pmatrix},\quad
        \begin{pmatrix}a&b\\0&0\end{pmatrix},\quad
        \begin{pmatrix}a&b\\c&0\end{pmatrix},\quad
        \begin{pmatrix}a&b\\c&d\end{pmatrix},
    $$
    with $a,b,c,d\in\F_2$. Their dimensions are $0,1,2,2,3,4$. By \Cref{prop:qmatroid-rectangles}, $\mM_{\mC_P}$ is a $q$-matroid exactly for $P\in\{\emptyset,(2),(2,2)\}$.

    Let $L=(\rho,|\cdot|,\YL(P_n))$ be the $\ZZ$-latroid with $\rho=(0,1,2,1,2,2)$. Its cyclic flats are $\emptyset$, $(2)$ and $(2,2)$. They form a chain, so $L$ is Schubert. The corresponding $q$-matroids are $\mU_{0,2}(2)$ for $P=\emptyset$, $\mM_1$ for $P=(2)$, and $\mU_{2,2}(2)$ for $P=(2,2)$. Here $\mM_1$ has rank function $V\mapsto\dim(V)-\dim(V\cap\langle(0,1)\rangle)$, and its space of loops is $\langle(0,1)\rangle$. Let $\mT:=\{\mU_{0,2}(2),\mM_1,\mU_{2,2}(2)\}$.

    The only flat of $\mU_{0,2}(2)$ is $\F_2^2$. The flats of $\mM_1$ are $\langle(0,1)\rangle$ and $\F_2^2$. Every subspace of $\F_2^2$ is a flat of $\mU_{2,2}(2)$. So the flats of each of these $q$-matroids are flats of the next one, and $\mT$ is a $q$-quotient collection.
\end{example}

\begin{corollary}\label{cor:cyclic-flats-qquotient}
    Let $L=(\rho,|\cdot|,\YL(P_n))$ be a $\ZZ$-latroid, and suppose that $\mM_{\mC_Z}$ is a $q$-matroid for every cyclic flat $Z\in\mZ_L$. Then $L$ is Schubert, and the $q$-matroids $\mM_{\mC_Z}$, $Z\in\mZ_L$, form a $q$-quotient collection.
\end{corollary}

\begin{proof}
    By \Cref{prop:qmatroid-rectangles}, every $Z\in\mZ_L$ is of the form $(h^s)$. Any two such diagrams are comparable, so $\mZ_L$ is a chain and $L$ is Schubert. If $Z\leq Z'$ in $\mZ_L$, then $\mM_{\mC_Z}$ is a $q$-quotient of $\mM_{\mC_{Z'}}$ by \Cref{prop:rank-difference}. Hence these $q$-matroids form a $q$-quotient collection.
\end{proof}

The hypothesis of \Cref{cor:cyclic-flats-qquotient} is restrictive: every cyclic flat must be a full-width rectangle $(h^s)$, and the $q$-matroids $\mM_{\mC_Z}$ then come from the nested codes in \Cref{rem:qdelta}\,(1) below. In particular, since $P_n$ is always a flat, it implies that $P_n$ is not cyclic, unless $P_n$ is of the form $(h^s)$.

\begin{remark}\label{rem:qdelta}
    We now look at \Cref{cor:cyclic-flats-qquotient} from the point of view of \cite{CJ25qanadeltamats}. Let $0\leq s\leq s'\leq t$, so that $\mM_s$ is a $q$-quotient of $\mM_{s'}$.
    \begin{enumerate}
        \item[(1)] Fix a basis of $\F_{q^h}$ over $\F_q$, and identify each vector $x\in\F_{q^h}^t$ with the $t\times h$ matrix over $\F_q$ whose $i$-th row is the coordinate vector of $x_i$. Then $\mC_{(h^s)}$ is identified with the $\F_{q^h}$-linear rank-metric code $\{x\in\F_{q^h}^t\mid x_{s+1}=\cdots=x_t=0\}$, which has generator matrix $(I_s\mid 0)$. Moreover, $\mM_s$ is the $q$-matroid of this code in the sense of \cite[Section~5.1]{CJ25qanadeltamats}. Indeed, if $A$ is a matrix whose columns form a basis of $V\leq E$, then the rank of $(I_s\mid 0)A$ over $\F_{q^h}$, which equals its rank over $\F_q$, is $\dim(V)-\dim(V\cap U_s^\perp)=r_s(V)$. So the cyclic flats in \Cref{cor:cyclic-flats-qquotient} give a chain of nested $\F_{q^h}$-linear codes, and the fact that $\mM_s$ is a $q$-quotient of $\mM_{s'}$ also follows from \cite[Proposition~5.2]{CJ25qanadeltamats}.
        \item[(2)] For a $q$-matroid on $E$ with rank function $r$, a subspace $V\leq E$ is \textbf{independent} if $r(V)=\dim(V)$, and \textbf{spanning} if $r(V)=r(E)$; see \cite[Definition~2.10]{CJ25qanadeltamats}. A \textbf{$q$-$\Delta$-matroid} on $E$ is a nonempty family of subspaces of $E$, called \textbf{feasible spaces}, that satisfies two exchange axioms; see \cite[Definition~3.1]{CJ25qanadeltamats}. If $\mathrm{Id}\colon\mN_1\to\mN_2$ is a strong map of $q$-matroids on $E$, then the subspaces that are independent in $\mN_1$ and spanning in $\mN_2$ are the feasible spaces of a \textbf{$q$-g-matroid} \cite[Definition~4.3]{CJ25qanadeltamats}, and every $q$-g-matroid is a $q$-$\Delta$-matroid; see \cite[Theorem~4.12]{CJ25qanadeltamats}. Since $\mM_s$ is a $q$-quotient of $\mM_{s'}$, the identity map $\mM_{s'}\to\mM_s$ is a strong map. Hence the subspaces of $E$ that are independent in $\mM_{s'}$ and spanning in $\mM_s$ are the feasible spaces of a $q$-g-matroid, and hence of a $q$-$\Delta$-matroid. We can describe them explicitly. A subspace $V\leq E$ is independent in $\mM_{s'}$ if and only if $V\cap U_{s'}^\perp=\{0\}$. Since $r_s(E)=s$ and, by Grassmann's identity, $r_s(V)=\dim(V+U_s^\perp)-(t-s)$, the subspace $V$ is spanning in $\mM_s$ if and only if $V+U_s^\perp=E$. Hence the feasible spaces are
        $$
            \{V\leq E\mid V\cap U_{s'}^\perp=\{0\}\text{ and }V+U_s^\perp=E\}.
        $$
        For instance, $U_k$ is feasible for every $s\leq k\leq s'$. This is a $q$-analogue of the family of sets $X$ with $\{1,\ldots,s\}\subseteq X\subseteq\{1,\ldots,s'\}$. This $q$-$\Delta$-matroid comes from a pair of nested $\F_{q^h}$-linear codes, so it is representable in the strong sense proposed in \cite[Section~5.1]{CJ25qanadeltamats}.
    \end{enumerate}
\end{remark}

\medskip

\section*{Acknowledgements}
G. N. Alfarano is supported by the Agence Nationale de la Recherche through grant number ANR-24-CPJ1-0075-01, and by Rennes Métropole through an AIS grant. S. Degen is supported by the Deutsche Forschungsgemeinschaft (DFG, German Research Foundation) through grants SFB-TRR 358/1 2023 -- 491392403 and SPP 2458 -- 539866293. S. Degen would like to thank Lukas K\"uhne for his supervision and guidance throughout the research process.

\medskip

\section*{Statement on the use of AI tools}
The authors used ChatGPT-6.1 Sol to assist with literature searches and language editing. AI tools were not used to develop mathematical proofs.

\bigskip

\bibliographystyle{abbrv}
\bibliography{references}
\clearpage

\appendix

\section{Enumeration of $\ZZ$-latroids on Young's lattice}\label{sec:enumeration}

In this appendix we enumerate all $\ZZ$-latroids $L=(\rho,|\cdot|,\YL(P_n))$ for any partition $P_n$ with $n\leq 4$. We denote any rank function $\rho$ on $\YL(P_n)$ by its tuple of rank values of the corresponding partitions in $\YL(P_n)$, which are in a prescribed total order, denoted by $\prec$. This total order is not to be confused with the usual partial order on $\YL(P_n)$, given by inclusion $\subseteq$.

\section*{Partitions of size $n=0$}

In this case we only have one option for $P_n$, namely $P_n=\emptyset$ with $\YL(P_n)=\{\emptyset\}$. Therefore, the only possible rank function is $(0)$. This rank function represents the uniform $\ZZ$-latroid of rank $0$.

\section*{Partitions of size $n=1$}

In this case we only have one option for $P_n$, namely $P_n=\ydiagram{1}$ with $\YL(P_n)=\{\emptyset \prec\ydiagram{1}\}$. Therefore, we have two possible rank functions, which are given by $(0,0)$ and $(0,1)$. This is due to the fact that by (L4) the rank can increase at most by one and by (L1) the rank of $\emptyset$ is always zero. Further note that, in this case, the axiom (L5) is superfluous, since $\YL(P_n)$ is already totally ordered with respect to inclusion. These two rank functions represent the uniform $\ZZ$-latroids of rank $0$ and $1$, respectively.

\section*{Partitions of size $n=2$}

In this case we have two options for $P_n$:

\begin{itemize}\setlength\itemsep{0.3cm}
    \item $P_n= \ydiagram{2}$ with lattice
    \[
        \YL(P_n)=\{\emptyset \prec\ydiagram{1}\prec\ydiagram{2}\}.
    \]
    This yields four possible rank functions:
    $$(0,0,0),\; (0,1,1),\; (0,1,2),\; (0,0,1).$$
    This is due to the fact that by (L4) the rank can increase at most by one and, moreover, that the axiom (L5) is superfluous, since $\YL(P_n)$ is already totally ordered with respect to inclusion.
    The first three represent the uniform $\ZZ$-latroids of rank $0,1,2$, respectively.

    \item $P_n= \ydiagram{1,1}$ with lattice
    \[
        \YL(P_n)=\{\emptyset \prec\ydiagram{1}\prec\ydiagram{1,1}\}.
    \]
    This yields exactly the same rank functions as $P_n=\ydiagram{2}$, by the same reasoning as in the former case.
\end{itemize}

\begin{remark}
    On $\YL(\ydiagram{2})$, the rank functions $(0,0,0)$ and $(0,0,1)$ define distinct $\ZZ$-latroids with the same independent partitions, namely
    only $\emptyset$, and the same circuits, namely only $\ydiagram{1}$.
    Hence neither the independent partitions nor the circuits
    determine a $\ZZ$-latroid on $\YL(P_n)$ in general.
    Their collections of flats, however, are different:
    $\{\ydiagram{2}\}$ and $\{\ydiagram{1},\ydiagram{2}\}$, respectively.
\end{remark}

\section*{Partitions of size $n=3$}

In this case we have three options for $P_n$:

\begin{itemize}\setlength\itemsep{0.3cm}
    \item $P_n= \ydiagram{3}$ with lattice
    \[
        \YL(P_n)=\{\emptyset \prec\ydiagram{1}\prec\ydiagram{2}\prec\ydiagram{3}\}.
    \]
    This yields eight possible rank functions:
    \[
        (0,0,0,0),\;(0,1,1,1),\;(0,1,2,2),\;(0,1,2,3),\;(0,0,0,1),\;(0,0,1,1),\;(0,0,1,2),\;(0,1,1,2).
    \]
    This is due to the fact that by (L4) the rank can increase at most by one and, moreover, that the axiom (L5) is superfluous, since $\YL(P_n)$ is already totally ordered with respect to inclusion.
    The first four represent the uniform $\ZZ$-latroids of rank $0,1,2,3$, respectively.

    \item $P_n= \ydiagram{1,1,1}$ with lattice
    \[
        \YL(P_n)=\{\emptyset \prec\ydiagram{1}\prec\ydiagram{1,1}\prec\ydiagram{1,1,1}\}.
    \]
    This yields exactly the same rank functions as $P_n=\ydiagram{3}$, by the same reasoning as in the former case.

    \item $P_n= \ydiagram{2,1}$ with lattice
    \[
        \YL(P_n)=\{\emptyset \prec\ydiagram{1}\prec\ydiagram{1,1}\prec\ydiagram{2}\prec\ydiagram{2,1}\}.
    \]
    Note that $\YL(P_n)$ is not totally ordered with respect to inclusion. Let us denote the above partitions by $Q_1,\ldots,Q_5$ from left to right in order. Furthermore, we denote by $(v_{Q_1},\ldots,v_{Q_5})$ their tuple of rank values. Then axioms (L1), (L4) and (L5) yield the following non-trivial inequalities for the values $v_{Q_1},\ldots,v_{Q_5}$:
    \begin{align*}
        &v_{Q_1}=0;\; 0\leq v_{Q_2}\leq1;\; v_{Q_3}\leq2;\; v_{Q_4}\leq2;\; v_{Q_5}\leq3\\
        &v_{Q_2}\leq v_{Q_3};\; v_{Q_2}\leq v_{Q_4};\;v_{Q_3},v_{Q_4}\leq v_{Q_5}\\
        &v_{Q_3}-v_{Q_2}\leq1;\; v_{Q_4}-v_{Q_2}\leq1;\; v_{Q_5}-v_{Q_3}\leq1;\; v_{Q_5}-v_{Q_4}\leq1\\
        &v_{Q_2}+v_{Q_5}\leq v_{Q_3}+v_{Q_4}.
    \end{align*}
    Here the third line comes from (L4) for the covers of $\YL(P_n)$, and the last line is (L5) for the only incomparable pair $Q_3,Q_4$, whose meet is $Q_2$ and whose join is $Q_5$.
    Now we distinguish between the two cases of $v_{Q_2}=0$ and $v_{Q_2}=1$.
    \begin{description}
        \item[Case $v_{Q_2}=0$]
        Then axiom (L4) implies that $0\leq v_{Q_3},v_{Q_4}\leq 1$ and hence also that $0\leq v_{Q_5}\leq 2$, which simplifies the inequality system as follows:
        \begin{align*}
            &v_{Q_1}=v_{Q_2}=0;\; 0\leq v_{Q_3},v_{Q_4}\leq1;\; v_{Q_5}\leq2\\
            &v_{Q_3},v_{Q_4}\leq v_{Q_5}\\
            &v_{Q_5}\leq v_{Q_3}+v_{Q_4}.
        \end{align*}
        Computing the integral solutions of the above inequality system yields:
        \[
            (0,0,0,0,0),\;(0,0,1,0,1),\;(0,0,0,1,1),\;(0,0,1,1,1),\;(0,0,1,1,2).
        \]

        \item[Case $v_{Q_2}=1$]
        Then axiom (L4) implies that $1\leq v_{Q_3},v_{Q_4}\leq 2$ and hence also that $1\leq v_{Q_5}\leq 3$, which simplifies the inequality system as follows:
        \begin{align*}
            &v_{Q_1}=0;\;v_{Q_2}=1;\; 1\leq v_{Q_3},v_{Q_4}\leq2;\; 1\leq v_{Q_5}\leq3\\
            &v_{Q_3},v_{Q_4}\leq v_{Q_5}\\
            &v_{Q_5}\leq v_{Q_3}+v_{Q_4}-1.
        \end{align*}
         Computing the integral solutions of the above inequality system yields:
        \[
            (0,1,1,1,1),\;(0,1,1,2,2),\;(0,1,2,1,2),\;(0,1,2,2,2),\;(0,1,2,2,3).
        \]
    \end{description}
    In both cases, the conditions $v_{Q_5}-v_{Q_3}\leq1$ and $v_{Q_5}-v_{Q_4}\leq1$ are not listed, since they follow from the last inequality together with $v_{Q_3},v_{Q_4}\leq v_{Q_2}+1$.
    Hence, we have the following ten rank functions for $P_n=\ydiagram{2,1}$:
    \begin{align*}
        &(0,0,0,0,0),\;(0,1,1,1,1),\;(0,1,2,2,2),\;(0,1,2,2,3),\\
        &(0,0,1,0,1),\;(0,0,0,1,1),\;(0,0,1,1,1),\;(0,0,1,1,2),\;(0,1,1,2,2),\;(0,1,2,1,2).
    \end{align*}
    The first four represent the uniform $\ZZ$-latroids of rank $0,1,2,3$, respectively.

\end{itemize}

\begin{remark}
    Another example of the aforementioned phenomenon is given by the rank functions $(0,1,1,1)$ and $(0,1,1,2)$ on $\YL(\ydiagram{3})$, which have the same independent partitions, namely $\emptyset$ and $\ydiagram{1}$, and the same circuits, namely only $\ydiagram{2}$. As before, their collections of flats are different. Indeed, they are $\{\emptyset,\ydiagram{3}\}$ and $\{\emptyset,\ydiagram{2},\ydiagram{3}\}$, respectively; see also \Cref{ex:chain-nongeometric}.
\end{remark}

\section*{Partitions of size $n=4$}

In this case we have five options for $P_n$ in total. The first two are those partitions such that the principal interval $ \YL(P_n)$ is given by a chain:

\begin{itemize}\setlength\itemsep{0.3cm}
    \item $P_n= \ydiagram{4}$ with lattice
    \[
        \YL(P_n)=\{\emptyset \prec\ydiagram{1}\prec\ydiagram{2}\prec\ydiagram{3}\prec\ydiagram{4}\}.
    \]
    This yields $16$ possible rank functions:
    \begin{align*}
        &(0,0,0,0,0),\;(0,1,1,1,1),\;(0,1,2,2,2),\;(0,1,2,3,3),\;(0,1,2,3,4);\;(0,0,0,0,1)\\&(0,0,0,1,1),\;(0,0,0,1,2),\;(0,0,1,1,1),\;(0,0,1,1,2),\;(0,0,1,2,2),\\
        &(0,0,1,2,3),\;(0,1,2,2,3),\;(0,1,1,2,2),\;(0,1,1,2,3),\;(0,1,1,1,2)
    \end{align*}

    This is due to the fact that by (L4) the rank can increase at most by one and, moreover, that the axiom (L5) is superfluous, since $\YL(P_n)$ is already totally ordered with respect to inclusion.
    The first five represent the uniform $\ZZ$-latroids of rank $0,1,2,3,4$, respectively.

    \item $P_n= \ydiagram{1,1,1,1}$ with lattice
    \[
        \YL(P_n)=\{\emptyset \prec\ydiagram{1}\prec\ydiagram{1,1}\prec\ydiagram{1,1,1}\prec\ydiagram{1,1,1,1}\}.
    \]
    This yields exactly the same rank functions as $P_n=\ydiagram{4}$, by the same reasoning as in the former case.
\end{itemize}

\begin{remark}
    Note that more generally for any $P_n$ such that $\YL(P_n)$ is a chain of $n+1$ partitions, we obtain $2^n$ many $\ZZ$-latroid rank functions.
\end{remark}

The next possibility is $P_n=\ydiagram{2,2}$ with lattice
    \[
        \YL(P_n)=\{\emptyset \prec\ydiagram{1}\prec\ydiagram{1,1}\prec\ydiagram{2}\prec\ydiagram{2,1}\prec\ydiagram{2,2}\}.
    \]
First note that $\YL(P_n)$ can be viewed as $\YL(\ydiagram{2,1})$ extended by a new top element, namely $\ydiagram{2,2}$. Hence, it is not hard to see that every rank function $\rho=(\rho_\emptyset,\rho_{(1)},\rho_{(1,1)},\rho_{(2)},\rho_{(2,1)},\rho_{(2,2)})$ in this case arises from rank functions $\rho^\prime=(\rho^\prime_\emptyset,\rho^\prime_{(1)},\rho^\prime_{(1,1)},\rho^\prime_{(2)},\rho^\prime_{(2,1)})$ of $\ydiagram{2,1}$ by choosing $\rho_{(2,2)}$ either to be $\rho_{(2,2)}=\rho^\prime_{(2,1)}$ or $\rho_{(2,2)}=\rho^\prime_{(2,1)}+1$. These are the only possible options for every choice of $\rho^\prime$, due to axiom (L4). Therefore, we have $2\cdot10=20$ different rank functions in total.\\

The last two options are $P_n=\ydiagram{3,1}$ and $P_n=\ydiagram{2,1,1}$, which yield the two isomorphic lattices:
\begin{align*}
   &\YL(P_n)=\{\emptyset \prec\ydiagram{1}\prec\ydiagram{1,1}\prec\ydiagram{2}\prec\ydiagram{3}\prec\ydiagram{2,1}\prec\ydiagram{3,1}\},\\
   &\YL(P_n)=\{\emptyset \prec\ydiagram{1}\prec\ydiagram{1,1}\prec\ydiagram{2}\prec\ydiagram{1,1,1}\prec\ydiagram{2,1}\prec\ydiagram{2,1,1}\},
\end{align*}
respectively. Transposition gives a lattice isomorphism between these two intervals. With the displayed orders, the corresponding rank tuples are obtained by interchanging their third and fourth entries. It therefore suffices to consider $P_n=(3,1)$.

Let us denote the above partitions by $Q_1,\ldots,Q_7$ from left to right in order. Furthermore, we denote by $(v_{Q_1},\ldots,v_{Q_7})$ their tuple of rank values. Then axioms (L1), (L4) and (L5) yield the following non-trivial inequalities for the values $v_{Q_1},\ldots,v_{Q_7}$:\vspace{0.3cm}

 \begin{center}
    \begin{tabular}{c|c|c|c}
          \textbf{Axioms (L1), (L4)} & \textbf{Axiom (L4)} & \textbf{Axiom (L4)}& \textbf{Axiom (L5)}\\ \hline
         $v_{Q_1}=0$ & $v_{Q_2}\leq v_{Q_3}$ & $v_{Q_3}-v_{Q_2}\leq 1$ & $v_{Q_2}+v_{Q_6}\leq v_{Q_3}+v_{Q_4}$\\
         $0\leq v_{Q_2}\leq1$ & $v_{Q_2}\leq v_{Q_4}$  & $v_{Q_4}-v_{Q_2}\leq 1$ & $v_{Q_2}+v_{Q_7}\leq v_{Q_3}+v_{Q_5}$ \\
         $0\leq v_{Q_3}\leq2$ & $v_{Q_3}\leq v_{Q_6}$ & $v_{Q_6}-v_{Q_3}\leq 1$ & $v_{Q_4}+v_{Q_7}\leq v_{Q_5}+v_{Q_6}$ \\
         $0\leq v_{Q_4}\leq2$ & $v_{Q_4}\leq v_{Q_5}$ & $v_{Q_5}-v_{Q_4}\leq 1$ &\\
         $0\leq v_{Q_5}\leq3$ & $v_{Q_4}\leq v_{Q_6}$ & $v_{Q_6}-v_{Q_4}\leq 1$ &\\
         $0\leq v_{Q_6}\leq3$ & $v_{Q_5}\leq v_{Q_7}$ & $v_{Q_7}-v_{Q_5}\leq 1$ &\\
         $0\leq v_{Q_7}\leq4$ & $v_{Q_6}\leq v_{Q_7}$ & $v_{Q_7}-v_{Q_6}\leq 1$ &
    \end{tabular}
\end{center}
\vspace{0.3cm}

Now we distinguish between the five cases $(v_{Q_2},v_{Q_3},v_{Q_4})=(0,0,0)$, $(v_{Q_2},v_{Q_3},v_{Q_4})=(0,1,0)$, $(v_{Q_2},v_{Q_3},v_{Q_4})=(0,0,1)$, $(v_{Q_2},v_{Q_3},v_{Q_4})=(0,1,1)$ and $v_{Q_2}=1$.
\begin{description}
        \item[Case $(v_{Q_2},v_{Q_3},v_{Q_4})=(0,0,0)$]
        Then the above inequality system simplifies as follows:
        \begin{align*}
            &v_{Q_1}=v_{Q_2}=v_{Q_3}=v_{Q_4}=v_{Q_6}=0;\; 0\leq v_{Q_5}\leq1\\
            &v_{Q_5}=v_{Q_7}.
        \end{align*}

        Computing the integral solutions of the above inequality system yields:
        \[
            (0,0,0,0,0,0,0),\;(0,0,0,0,1,0,1).
        \]

        \item[Case $(v_{Q_2},v_{Q_3},v_{Q_4})=(0,1,0)$]
        Then the above inequality system simplifies as follows:
        \begin{align*}
            &v_{Q_1}=v_{Q_2}=v_{Q_4}=0;\; v_{Q_3}=v_{Q_6}=1,\; 0\leq v_{Q_5}\leq1;\; 1\leq v_{Q_7}\leq2\\
            &v_{Q_5}\leq v_{Q_7}\\
            &v_{Q_7}\leq v_{Q_5}+1.
        \end{align*}

        Computing the integral solutions of the above inequality system yields:
        \[
            (0,0,1,0,0,1,1),\;(0,0,1,0,1,1,1),\;(0,0,1,0,1,1,2).
        \]

        \item[Case $(v_{Q_2},v_{Q_3},v_{Q_4})=(0,0,1)$]
        Then the above inequality system simplifies as follows:
        \begin{align*}
            &v_{Q_1}=v_{Q_2}=v_{Q_3}=0;\; v_{Q_4}=v_{Q_6}=1,\; 1\leq v_{Q_5}\leq2;\; 1\leq v_{Q_7}\leq2\\
            &v_{Q_5}=v_{Q_7}.
        \end{align*}

        Computing the integral solutions of the above inequality system yields:
        \[
            (0,0,0,1,1,1,1),\;(0,0,0,1,2,1,2).
        \]

        \item[Case $(v_{Q_2},v_{Q_3},v_{Q_4})=(0,1,1)$]
        Then the above inequality system simplifies as follows:
        \begin{align*}
            &v_{Q_1}=v_{Q_2}=0;\; v_{Q_3}=v_{Q_4}=1,\; 1\leq v_{Q_5},v_{Q_6}\leq2;\; 1\leq v_{Q_7}\leq3\\
            &v_{Q_5},v_{Q_6}\leq v_{Q_7}\\
            &v_{Q_7}\leq v_{Q_5}+1;\; v_{Q_7}\leq v_{Q_5}+v_{Q_6}-1.
        \end{align*}

        Computing the integral solutions of the above inequality system yields:
        \[
            (0,0,1,1,1,1,1),\;(0,0,1,1,1,2,2),\;(0,0,1,1,2,1,2),\;(0,0,1,1,2,2,2),\;(0,0,1,1,2,2,3).
        \]

        \item[Case $v_{Q_2}=1$]
        For each of the previous four cases we can shift their integral solutions by the vector $(0,1,1,1,1,1,1)$ to obtain all possible integral solutions in this situation.
    \end{description}
    Since all of these integral solutions are clearly distinct, we therefore have $2\cdot(2+3+2+5)=24$ different rank functions in total.

\begin{remark}
We would like to emphasize that the individual inequality systems used to compute all of these $\ZZ$-latroid rank functions on the different lattices $\YL(P_n)$ determine constraints of analogues of the $q$-rank polytopes described in \cite{alfaranodegen2025rankpolytope}. These analogue polytopes form the sets of all respective $\R$-latroids $(\rho,|\cdot|,\YL(P_n))$, and hence, the respective integral solutions correspond to their respective lattice points.
\end{remark}

\end{document}